\documentclass[12pt]{amsart}

\usepackage{amsmath, amssymb, amsfonts, amsthm, amscd}
\usepackage{manfnt}
\usepackage{mathtools}
\usepackage{fullpage}
\usepackage[all]{xy}
   \SelectTips{cm}{10}
\usepackage{txfonts}
\usepackage{enumitem}
\newtheorem{thm}{Theorem}[section]
\newtheorem{lemma}[thm]{Lemma}
\newtheorem{prop}[thm]{Proposition}
\newtheorem{cor}[thm]{Corollary}

\newtheorem{conj}[thm]{Conjecture}

\newtheorem{prop-conj}[thm]{Proposition-Conjecture}

\theoremstyle{definition}
\newtheorem{defn}[thm]{Definition}

\theoremstyle{remark}
\newtheorem{rmk}[thm]{Remark}

\theoremstyle{remark}
\newtheorem{notation}[thm]{Notation}

\theoremstyle{remark}

\theoremstyle{remark}

\newcommand{\Q}{\mathbb{Q}}
\newcommand{\Qb}{\overline{\mathbb{Q}}}
\newcommand{\Z}{\mathbb{Z}}
\newcommand{\CC}{\mathbb{C}}
\newcommand{\RR}{\mathbb{R}}

\newcommand{\Ql}{\mathbb{Q}_{\ell}}
\newcommand{\Qlb}{\overline{\mathbb{Q}}_\ell}

\DeclareMathOperator{\Hom}{Hom}

\DeclareMathOperator{\End}{End}
\DeclareMathOperator{\Aut}{Aut}

\DeclareMathOperator{\Ad}{Ad}

\DeclareMathOperator{\im}{im}

\DeclareMathOperator{\Rep}{Rep}
\DeclareMathOperator{\Spec}{Spec}
\DeclareMathOperator{\id}{id}
\DeclareMathOperator{\Bun}{Bun}
\DeclareMathOperator{\Gr}{Gr}
\DeclareMathOperator{\Fl}{Fl}

\DeclareMathOperator{\Sat}{Sat}
\DeclareMathOperator{\IC}{IC}
\DeclareMathOperator{\Loc}{Loc}
\DeclareMathOperator{\Stab}{Stab}

\DeclareMathOperator{\Vect}{Vect}
\DeclareMathOperator{\ev}{ev}
\DeclareMathOperator{\IH}{IH}
\DeclareMathOperator{\Pic}{Pic}
\DeclareMathOperator{\Frac}{Frac}

\newcommand{\gal}[1]{\Gamma_{#1}} 
\newcommand{\Gal}{\mathrm{Gal}} 

\newcommand{\into}{\hookrightarrow}
\newcommand{\onto}{\twoheadrightarrow}

\newcommand{\mc}{\mathcal}
\newcommand{\mf}{\mathfrak}
\newcommand{\mr}{\mathrm}
\newcommand{\mbf}{\mathbf}

\newcommand{\tG}{\widetilde{G}}
\newcommand{\tT}{\widetilde{T}}
\newcommand{\tZ}{\widetilde{Z}}
\newcommand{\tB}{\widetilde{B}}
\newcommand{\tH}{\widetilde{H}}

\newcommand{\Ring}{\mathrm{Ring}}

\newcommand{\tPone}{\widetilde{\mathbf{P}}_0(1)}
\newcommand{\tI}{\widetilde{\mathbf{I}}}
\newcommand{\tPinf}{\widetilde{\mathbf{P}}_{\infty}}
\newcommand{\tBun}{\widetilde{\mathrm{Bun}}}
\newcommand{\tHk}{\widetilde{\mathrm{Hk}}}
\newcommand{\lefth}{\overleftarrow{h}}
\newcommand{\righth}{\overrightarrow{h}}
\newcommand{\tGR}{\widetilde{\mathrm{GR}}}
\newcommand{\tE}{\widetilde{\mathcal{E}}}
\newcommand{\tmfG}{\widetilde{\mathfrak{G}}^U_{\leq \lambda}}
\newcommand{\twoK}{{}^{(2)}K}
\newcommand{\twoP}{{}^{(2)}\mathbf{P}_0}
\newcommand{\twoZ}{{}^{(2)}Z_G}
\newcommand{\twoA}{{}^{(2)}\mathrm{A}}

\newcommand{\bs}{\backslash}

\begin{document}
\title{Generalized Kuga-Satake theory and rigid local systems, II: rigid Hecke eigensheaves}
\thanks{This paper could not exist without the beautiful work \cite{yun:exceptional} of Zhiwei Yun. I am moreover very grateful to Yun for generously sharing his understanding of \cite{yun:exceptional} with me. I thank Bhargav Bhatt for answering a question about stratifications, and Mark de Cataldo for his feedback on a draft of this paper. I thank the Institut Math\'{e}matique de Jussieu for its hospitality during part of the preparation of this paper. This work was supported by NSF grant DMS-1303928.}
\author{Stefan Patrikis}
\email{patrikis@math.mit.edu}
\address{MIT Department of Mathematics \\ Building E18\\ 77 Massachusetts Avenue \\Cambridge, MA 02139}
\date{May 2014}
\begin{abstract}
This paper uses rigid Hecke eigensheaves, building on Yun's work on the construction of motives with exceptional Galois groups, to produce the first robust examples of `generalized Kuga-Satake theory' outside the Tannakian category of motives generated by abelian varieties. To strengthen our description of the `motivic' nature of Kuga-Satake lifts, we digress to establish a result that should be of independent interest: for any quasi-projective variety over a (finitely-generated) characteristic zero field, the associated weight-graded of its intersection cohomology arises from a motivated motive in the sense of Andr\'{e}, and in particular from a classical homological motive if one assumes the Standard Conjectures. This extends work of de Cataldo and Migliorini.
\end{abstract}
\maketitle

\section{Background: generalized Kuga-Satake theory}
The aim of this paper is to produce non-trivial examples of the generalized Kuga-Satake theory proposed in \cite{stp:variationsarxiv}. The classical Kuga-Satake construction is a miracle of Hodge theory that associates to any complex $K3$ surface $X$ a complex abelian variety $\mr{KS}(X)$ and an inclusion of $\Q$-Hodge structures
\[
H^2(X, \Q) \subset H^1(\mr{KS}(X), \Q)^{\otimes 2}.
\]
This construction takes its clearest conceptual form within the motivic Galois formalism. Let $\mc{M}_{\CC}^{\mr{hom}}$ denote the category of pure motives over $\CC$ for homological equivalence. Assuming the Standard Conjectures, this is a neutral Tannakian category over $\Q$ with fiber functor given by Betti cohomology:
\[
H_B \colon \mc{M}_{\CC}^{\mr{hom}} \to \Vect_{\Q}.
\]
Let $\mc{G}^{\mr{hom}}_{\CC}= \Aut^{\otimes}(H_B)$ denote the corresponding Tannakian group. Then we can phrase the Kuga-Satake construction as follows: the motive $H^2(X)$ admits a (symmetric) polarization, hence (normalizing by a Tate-twist to weight zero) corresponds to a motivic Galois representation $\rho \colon \mc{G}_{\CC}^{\mr{hom}} \to \mr{SO}(H^2_B(X)(1))$.\footnote{That this motivic Galois representation should be special orthogonal, rather than merely orthogonal, is non-trivial; it follows from work of Andre: \cite{andre:hyperkaehler}.} The motive $H^1(\mr{KS}(X))$ then \textit{is} the motivic Galois representation corresponding to the composite $r \circ \tilde{\rho}$ in the diagram
\begin{equation}\label{classicalKS}
\xymatrix{
& \mr{GSpin}(H^2_B(X)(1)) \ar[r]^r \ar[d] & \mr{GL}(C^+(H^2_B(X)(1))) \\
\mc{G}_{\CC}^{\mr{hom}} \ar@{-->}[ur]^{\tilde{\rho}} \ar[r]_-{\rho} & \mr{SO}(H^2_B(X)(1)),
}
\end{equation}
where $\tilde{\rho}$ is a suitable lift of $\rho$, and $r$ is the natural representation of $\mr{GSpin}$ on the even Clifford algebra. The strongest possible version of the Kuga-Satake construction is the statement that such a lift $\tilde{\rho}$ exists; this is far from known at present, as it implicitly includes deep cases of the Lefschetz standard conjecture. A weaker, but still highly non-trivial, analogue is known when $\mc{G}_{\CC}^{\mr{hom}}$ is replaced by the motivic Galois group of Andr\'{e}'s category of motives for motivated cycles; see \cite{andre:hyperkaehler}.

But the formulation itself is highly suggestive, pointing towards deep and largely unexplored generalizations, some of whose essential difficulties are orthogonal to the usual impenetrable conjectures of algebraic and arithmetic geometry: Lefschetz, Hodge, Tate, etc. In what follows we will work with motives over number fields and their $\ell$-adic realizations, rather than motives over $\CC$ and their Hodge-Betti realizations, but there are analogues of the results of this paper in the latter setting. We now state a conjecture that captures the most refined form of a `generalized Kuga-Satake theory' for motives over number fields. For two number fields $F$ and $E$, we let $\mc{M}_{F, E}$ denote the category of motives for motivated cycles over $F$ with coefficients in $E$; it is (unconditionally) neutral Tannakian over $E$, and by choosing an embedding $F \into \CC$, the ($E$-linear) Betti fiber functor gives us its motivic Galois group $\mc{G}_{F, E}$ (see \cite{andre:motivated} for background).
\begin{conj}[See \S 4.3 of \cite{stp:variationsarxiv}]\label{KSconj}
Let $\tH \to H$ be a surjection of linear algebraic $E$-groups whose kernel is equal to a central torus in $\tH$, and let
\[
\rho \colon \mc{G}_{F, E} \to H
\]
be a motivic Galois representation. Then if either $F$ is totally imaginary, or the `Hodge numbers' of $\rho$ satisfy the (necessary) parity condition of \cite[Proposition 5.5]{stp:parities}, then there exists a finite extension $E'/E$ and a lifting of motivic Galois representations
\[
\xymatrix{
& \tH_{E'} \ar[d] \\
\mc{G}_{F, E'} \ar[r]_{\rho \otimes_E E'} \ar[ur]^{\tilde{\rho}} & H_{E'}.
}
\]
\end{conj}
For a leisurely overview of this conjecture, see the introduction to \cite{stp:KSRLS1:middleconv}; for a detailed discussion of the arithmetic evidence, see \cite{stp:variationsarxiv}. Even working with motivated rather than homological motives, this conjecture is highly refined: in the classical setting of diagram (\ref{classicalKS}), the existence of such a $\tilde{\rho}$ requires not only the existence of $\mr{KS}(X)$, but also the full force of the theorem of Deligne-Andr\'{e} that Hodge cycles on abelian varieties are motivated.\footnote{For $F=\CC$, that is. When $F$ is a number field, the existence of $\tilde{\rho}$ requires the Tate conjecture for abelian varieties, and rather delicate descent arguments: see \cite[\S 4.2]{stp:variationsarxiv}.} At first approximation, though, we can replace Conjecture \ref{KSconj} with the following variant:
\begin{defn}\label{genKSdefn}
Setting $\gal{F}= \Gal(\overline{F}/F)$ for an algebraic closure $\overline{F}$ of $F$, let $\rho \colon \gal{F} \to H(\Qlb)$ be a geometric Galois representation valued in an arbitrary linear algebraic group $H$ over $\Qlb$. 
\begin{itemize}
\item We say that $\rho$ is weakly motivic if there exists a \textit{faithful} representation $r \colon H \into \mr{GL}(V_r)$ such that $r \circ \rho$ is isomorphic to the $\iota \colon E \into \Qlb$ realization $H_{\iota}(M)$ of some object $M$ of $\mc{M}_{F, E}$. 
\item Suppose that we are given such a weakly motivic $\rho \colon \gal{F} \to H(\Qlb)$, and let $\tilde{\rho}$ be a geometric lift to $\tH$:
\begin{equation}\label{galoisside}
\xymatrix{
& \tH(\Qlb) \ar[d] \\
\gal{F} \ar[ur]^{\tilde{\rho}} \ar[r]_-{\rho} & H(\Qlb).
}
\end{equation}
(That such geometric lifts typically exist is \cite[Theorem 3.2.10]{stp:variationsarxiv} and \cite[Proposition 5.5]{stp:parities}.\footnote{This should be contrasted with the situation in which the kernel of $\tH \to H$ is an isogeny, where geometric lifts, even after allowing a finite base-change on $F$, need not exist: for a simple example, consider the case $\mr{SL}_2 \to \mr{PGL}_2$ in which $\rho$ is the projective representation associated to the Tate module of an elliptic curve (or even more simply, consider multiplication by $N >1$ on $\mathbb{G}_m$, and let $\rho$ be the cyclotomic character). For the full story, see \cite{wintenberger:relevement}.}) We say that $\tilde{\rho}$ satisfies the generalized Kuga-Satake property if it is weakly motivic as $\tH$-representation.
\end{itemize}
\end{defn}
In sum, our aim in establishing certain cases of this `generalized Kuga-Satake property' is to verify (motivated refinements of) certain cases of the Fontaine-Mazur conjecture.

With this framework in place, we can introduce the particular setting of this paper. Our aim is to study certain families of weakly motivic $\rho \colon \gal{F} \to H(\Qlb)$ for which it is possible to find lifts $\tilde{\rho} \colon \gal{F} \to \tH(\Qlb)$ satisfying the generalized Kuga-Satake property. Outside of the context of the classical Kuga-Satake construction, where $\rho$ is the representation on $H^2(X_{\overline{F}}, \Qlb)$, for $X/F$ a $K3$ surface--or closely related examples in which the motives in question are still generated by motives of abelian varieties\footnote{For an `axiomatic' generalization of this context, see \cite{andre:hyperkaehler}, which, for instance, further allows $X$ to be a hyperk\"{a}hler variety, or a cubic four-fold.}--there were no non-trivial examples of such lifting until the paper \cite{stp:KSRLS1:middleconv}. But that paper is restricted to low-dimensional examples in which $\tH= \mr{GSpin}_5 \to H= \mr{SO}_5$, and relies heavily on low-dimensional coincidences in the Dynkin classification. Thus the primary desiderata for our examples are that:
\begin{enumerate}[label= D.\arabic*, ref= D.\arabic*]
\item \label{nonabelian} the motives in question not lie in the Tannakian subcategory of $\mc{M}_F$ generated by abelian varieties and Artin motives;
\item \label{allranks} the examples exist in arbitrary rank, or at least for `interesting' groups $H$;
\item \label{nontriv} the lift $\tilde{\rho}$ should not be realizable within the Tannakian category of geometric representations generated by $\rho$, characters, and Artin representations.
\end{enumerate} 
We make explicit this last desideratum just to point out that for some choices of $\tH$, for instance $\tH= H \times \mathbb{G}_m$, the existence of a weakly motivic lift $\tilde{\rho}$ is completely trivial. Condition \ref{nontriv} is a way to ensure the results we prove have non-trivial content.

The examples of this paper meet all three criteria of interest. For our $\rho$ we take the remarkable weakly motivic Galois representations constructed by Yun in \cite[Theorem 4.2, Proposition 4.6]{yun:exceptional}. Let us recall a somewhat simplified version of the main result of \cite{yun:exceptional}. Let $G$ be a split, simple, simply-connected group of type $A_1$, $D_n$ with $n$ even, $G_2$, $E_7$, or $E_8$, and let $G^\vee$ denote the split $\Q$-form of its dual group. We have to say a word about the coefficients of the Galois representations and motives. For definiteness, fix an embedding $\Qb \into \Qlb$, implicit whenever we take `the' $\ell$-adic realization of a motive with coefficients in $\Qb$, and let $\imath$ be a square-root of $-1$ in $\Qb$. All the local systems considered can be arranged to have coefficients in the (possibly trivial) extension $\Ql'= \Ql(\imath)$. The motives will have coefficients in the subfield $\Q' \subset \Qb$ given by
\begin{equation}
\Q'= \begin{cases}
\text{$\Q$ in types $D_{4m}$, $G_2$, $E_8$; and}\\
\text{$\Q(\imath)$ in types $A_1$, $D_{4m+2}$, $E_7$.}
\end{cases}   
\end{equation}
There is a certain two-fold cover $\twoZ \onto Z_G$ (see Definition \ref{thefinitegroup} and Lemma \ref{grouptheory}) of the center $Z_G$ of $G$--regard $\twoZ$ as a group scheme over $\Q$--and we call a character 
\[
\chi \colon \twoZ(\Qb) \to \Qb^\times
\] 
\textit{odd} if it is non-trivial on the kernel of $\twoZ \to Z_G$.
\begin{thm}[\cite{yun:exceptional}]\label{yun}
For any odd character $\chi \colon \twoZ(\Qb) \to \Qb^\times$, there exists a local system
\[
\rho_{\chi} \colon \pi_1(\mathbb{P}^1_{\Q}-\{0,1,\infty\}) \to G^\vee(\Ql')
\]
whose geometric monodromy group is $G^\vee$, except in type $D_{2m}$, in which case the geometric monodromy group is $\mr{SO}_{4m-1}$. For all number fields $F$ such that
\[
F \supseteq
\begin{cases}
\text{$\Q$ if $G$ is type $D_{4m}$, $G_2$, or $E_8$;} \\
\text{$\Q(\sqrt{-1})$ if $G$ is type $A_1$, $D_{4m+2}$, or $E_7$,}
\end{cases}
\] 
and all specializations $t \colon \Spec F \to \mathbb{P}^1-\{0,1,\infty\}$, the pull-back $\rho_{\chi, t} \colon \gal{F} \to G^\vee(\Ql')$ is weakly motivic. To be precise, the composition of $\rho_{\chi, t}$ with the quasi-minuscule representation of $G^\vee$ is isomorphic to the $\Ql'$-realization of an object of $\mc{M}_{F, \Q'}$.
\end{thm}
We can now state the first main result of this paper. There is a minor technicality in the phrasing of this theorem that results very naturally from the way the geometric Satake isomorphism descends to number fields: see \S \ref{satake} for a careful explanation. Namely, for any connected reductive group $H$, let $\rho^\vee$ denote the usual half-sum of the positive coroots (for any choice of based root datum), and set $H_1= (H \times \mathbb{G}_m)/ \langle (2\rho^\vee(-1) \times -1) \rangle$. In the case $H=G^\vee$, to avoid cluttered notation we write $G^\vee_1$ for $(G^\vee)_1$; this should not cause any confusion. Yun's construction is most naturally viewed as the construction of a local system 
\[
\rho_{\chi} \colon \pi_1(\mathbb{P}^1_{\Q}-\{0,1,\infty\}) \to G^\vee_1(\Ql')= (G^\vee \times \mathbb{G}_m)(\Ql')
\]
in which the $G^\vee$-component is as in Theorem \ref{yun}, and the $\mathbb{G}_m$ component is the cyclotomic character; the equality here uses the fact that $G$ is simply-connected.
\begin{thm}\label{mainthm}
Let $\tH \onto G^\vee$ be any surjection of split connected reductive groups with kernel equal to a central torus in $\tH$. Then:
\begin{enumerate}
\item There exists a local system $\tilde{\rho}_{\chi} \colon \pi_1(\mathbb{P}^1_{\Q(\sqrt{-1})}-\{0,1,\infty\}) \to \tH_1(\Ql')$ lifting $\rho_{\chi}$, i.e. such that the diagram
\[
\xymatrix{
& \tH_1(\Ql') \ar[d] \\
\pi_1(\mathbb{P}^1_{\Q(\sqrt{-1})}-\{0,1,\infty\}) \ar[ur]^{\tilde{\rho}_{\chi}} \ar[r]_-{\rho_{\chi}} & G^\vee_1(\Ql')
}
\]
commutes. When $G$ is type $D_{4m}$, we may replace $\Q(\sqrt{-1})$ by $\Q$ in this assertion.
\item For all number field specializations $t \colon \Spec F \to \mathbb{P}^1-\{0,1,\infty\}$ (assuming $F \supset \Q(\sqrt{-1})$ in types $A_1$, $D_{4m+2}$, and $E_7$), $\tilde{\rho}_{\chi, t}$ is weakly motivic, i.e. satisfies the generalized Kuga-Satake property.
\end{enumerate}
\end{thm}
The real content of this result is for $G$ of types $D_{2m}$ and $E_7$; when $\pi_1(G)= \{1\}$ (types $G_2$, $E_8$), there can never be any generalized Kuga-Satake lift satisfying criterion \ref{nontriv}. In type $A_1$, the construction is not completely trivial, but the motives in question are generated by abelian varieties and Artin motives, so fail to satisfy our criterion \ref{nonabelian}.\footnote{Also, in this case, a more elementary construction of the lift can be achieved using Katz's theory (\cite{katz:rls}) of rigid local systems; this is a simple case of the strategy of \cite{stp:KSRLS1:middleconv}.} But in the essential cases of types $D_{2m}$ and $E_7$, all of our desiderata are met, the key point being that for suitable choice of $\tH$, the group $\tH_1$ has irreducible representations restricting to each of the minuscule representations of the simply-connected cover $G^\vee_{\mr{sc}}$ of $G^\vee$; these are representations not possessed by the original (adjoint) group $G^\vee$. See \S \ref{minuscule} for details.

We now briefly summarize the approach to constructing the lifted local systems $\tilde{\rho}_{\chi}$. Yun's $\rho_{\chi}$ is constructed as the eigen-local system associated to a Hecke eigensheaf on a certain moduli space $\Bun$ of $G$-bundles on $\mathbb{P}^1$ with level structure at the points $\{0,1,\infty\}$. Simply put, we enlarge the center of the semi-simple group $G$ to form a reductive group $\tG$ (whose dual group $\tG^\vee$ plays the role of $\tH$ above); then we study an analogous moduli space $\tBun$ of $\tG$-bundles with level structure, and show that Yun's eigensheaves can be extended to eigensheaves on $\tBun$. The weakly motivic nature of the lifts $\tilde{\rho}_{\chi, t}$ is realized in the (restricting to the interesting cases in type $A_1$, $D_{2m}$, $E_7$) minuscule representations of $\tG^\vee$ (or rather, of $\tG^\vee_1$); as in \cite{yun:exceptional}, the motives themselves are closely related to the (intersection) cohomology of certain open subvarieties of affine Schubert varieties.

To put this approach in perspective, let us note that it is a geometric analogue of the classical automorphic construction parallel to the lifting problem (\ref{galoisside}). Namely, extending an automorphic representation of $G$ to $\tG$ heuristically corresponds to lifting a representation $\mc{L}_F \to G^\vee(\CC)$ of the `automorphic Langlands group' $\mc{L}_F$ to $\tG^\vee(\CC)$. We are carrying out an analogue for certain Hecke eigensheaves, being careful to retain hold of the explicit `motivic' nature of the corresponding eigen-local systems.

In fact, we prove something considerably stronger than Theorem \ref{mainthm}, strengthening the `motivic' result even in Yun's original context. Rather than showing (as in part 2 of Theorem \ref{mainthm}) that the $\tilde{\rho}_{\chi, t}$ (or $\rho_{\chi, t}$) are weakly motivic, we show (Theorem \ref{motivated}) that for \textit{any} finite-dimensional representation $r$ of $\tH_1$, $r \circ \tilde{\rho}_{\chi, t}$ is motivated. The content of this assertion is the following: the arguments showing that $\rho_{\chi, t}$ and $\tilde{\rho}_{\chi, t}$ are weakly motivic rest on the fact that quasi-miniscule and miniscule affine Schubert varieties have very mild singularities (punctual in the quasi-miniscule case; none at all in the miniscule case). For such varieties (and their close cousins that appear in the proof), we can in quite elementary terms describe their intersection cohomology groups as the $\ell$-adic realizations of motivated motives. The claim that \textit{all} $r \circ \tilde{\rho}_{\chi, t}$ are motivated depends on a similar description, but for varieties with singularities as bad as those of any affine Schubert variety. This essentially means we need a `motivated' description of the intersection cohomology $\IH^*(Y_{\bar{k}}, \Ql)$ of an arbitrarily singular, and not necessarily projective, variety $Y$ over a characteristic zero field $k$; to be precise, since motivated motives do not reflect `mixed' behavior, we prove such an assertion for the associated weight-graded $\Gr^W_{\bullet} \IH^*(Y_{\bar{k}}, \Ql)$. This is deduced as a consequence of a stronger `motivated' variant of the decomposition theorem, and especially from a `motivated support decomposition': see Theorem \ref{motivateddecomposition} and Corollary \ref{splitting}. Here is the statement for intersection cohomology:
\begin{thm}[see Corollary \ref{ICmotpf}]\label{ICmotintro}
Let $k$ be a finitely-generated field of characteristic zero, and let $Y/k$ be any quasi-projective variety. Then there is an object $M \in \mc{M}_k$ whose $\ell$-adic realization is isomorphic as $\gal{k}$-representation to $\Gr^W_i \IH^m(Y_{\bar{k}}, \Ql)$. If $\Gamma$ is a finite group scheme over $k$ acting on $Y$, and $e \in \Qb[\Gamma(\bar{k})]^{\gal{k}}$ is an idempotent, then for any embedding $\Qb \into \Qlb$ there is an object of $\mc{M}_{k, \Qb}$ whose $\Qb \into \Qlb$-realization is isomorphic as $\gal{k}$-representation to $\Gr^W_i e(\IH^m(Y_{\bar{k}}, \Qlb)$. 

The same holds for intersection cohomology with compact supports.
\end{thm}
When $Y$ is projective, in which case $\IH^m(Y_{\bar{k}}, \Ql)$ is pure of weight $m$, and $k$ is algebraically closed, this\footnote{Not exactly, of course, since as we have phrased the result the theorem is vacuous for $k$ algebraically closed; in that case substitute for the $\ell$-adic cohomology the collection of Betti, de Rham, and $\ell$-adic realizations.} is a recent result of de Cataldo-Migliorini (\cite[Theorem 3.2.2]{decataldo-migliorini:decomproj}), part of their beautiful series of papers (see for instance \cite{decataldo-migliorini:htalgmaps}, \cite{decataldo-migliorini:pervfil}, \cite{decataldo:pervfil}) re-establishing the decomposition theorem and its associated mixed Hodge-theoretic package by `geometric,' rather than `sheaf-theoretic,' methods. These papers chart a fundamental advance in our understanding of the geometry of perverse sheaves, and I expect will find many more, and far deeper, motivic applications than the one here. Since the arguments establishing Theorem \ref{ICmotintro} are independent of those of the rest of this paper, I refer the reader to \S \ref{overview} for a fuller introduction, and for an overview of the approach to Theorem \ref{ICmotintro}. Also see Remark \ref{finalrmk} for additional applications, such as a $p$-adic de Rham comparison isomorphism for intersection cohomology.
\section{Bundles with level structure}\label{levelstructure}
In this section only, we allow $G$ to be any connected reductive group over a field $k$, and $X$ to be any smooth projective geometrically connected curve over $k$. Our aim is to review a construction from \cite[\S 4.2]{yun:globalspringer} of moduli spaces of $G$-bundles on $X$ with level structure at a finite set $S= \{x_1, \ldots, x_n\} \subset X(k)$ of $k$-points. Here and throughout, we denote by $LG$ and $L^+G$ the `abstract' loop group and positive loop group of $G$, i.e. the functors of $k$-algebras given by $R \mapsto G\left(R((t))\right)$ and $R \mapsto G(R[[t]])$ (a group ind-scheme and pro-algebraic group, respectively, over $k$), where $t$ is a formal parameter. Now let $x$ be a closed point of $X$, and denote by $\mc{O}_x$ the complete local ring of $X$ at $x$, with fraction field $\mc{K}_x$. Then we denote by $L_x G$ and $L^+_x G$ the functors $R \mapsto G \left(R \widehat{\otimes}_{\kappa(x)} \mc{K}_x \right)$ and $R \mapsto G \left(R \widehat{\otimes}_{\kappa(x)} \mc{O}_x \right)$.
\begin{defn}
Let $\Bun_{G, S, \infty} \to \Ring_k$ be the stack associated to the following prestack $\Bun_{G, S, \infty}^{pre}$ over $k$: for any $k$-algebra $R$, $\Bun_{G, S, \infty}^{pre}(R)$ is the groupoid of triples $(\alpha, \mc{P}, \tau)$ where
\begin{itemize}
\item $\alpha= (\alpha_{x_i})_{i=1, \ldots, n}$ is a collection of local coordinates $\alpha_{x_i} \colon R[[t]] \xrightarrow{\sim} \widehat{\mc{O}}_{x_i}$ (here we regard $x_i$ as an $R$-point $x_i \colon \Spec R \to X_R$ and take the formal completion of $X_R$ along the graph $\Gamma(x_i)$);
\item $\mc{P}$ is a $G$-torsor on $X_R$;
\item $\tau= (\tau_{x_i})_{i= 1, \ldots, n}$ is a collection of full level structures (abbreviating $\mc{D}_{x_i}= \Spec (\widehat{\mc{O}}_{x_i})$)
\[
\tau_{x_i} \colon G \times \mc{D}_{x_i} \xrightarrow{\sim} \mc{P}|_{\mc{D}_{x_i}}.
\]
\end{itemize}
\end{defn}
Let $\Aut_{\mc{O}}$ denote the pro-algebraic group of continuous automorphisms of $k[[t]]$. The semi-direct product
\[
(LG \rtimes \Aut_{\mc{O}})^n
\]
acts on the right on $\Bun_{G, S, \infty}$ as follows.
\begin{defn}\label{levelaction}
For $g= (g_i)_{i=1, \ldots n} \in G(R((t))^n$ and $\sigma= (\sigma_i)_{i=1, \ldots, n} \in \Aut(R[[t]])^n$, and $(\alpha, \mc{P}, \tau) \in \Bun_{G, S, \infty}^{pre}(R)$, let $(g, \sigma)$ act on $(\alpha, \mc{P}, \tau)$ by
\[
R_{g, \sigma}(\alpha, \mc{P}, \tau)= (\alpha \circ \sigma, \mc{P}^{g}, \tau^{g}),
\] 
where:
\begin{itemize}
\item $\alpha \circ \sigma= (\alpha_{x_i} \circ \sigma_i)_i$;
\item $\mc{P}^g$ is the $G$-bundle on $X_R$ obtained by gluing $\mc{P}|_{X_R- \cup_i \Gamma(x_i)}$ to the trivial $G$-bundles on the the completions $\mc{D}_{x_i}= \widehat{\mc{O}}_{x_i}$ along the punctured discs $\mc{D}_{x_i}^\times$ via the isomorphisms
\[
G \times \mc{D}_{x_i}^\times \xrightarrow{\alpha_{x_i} \circ g_i \circ \alpha_{x_i}^{-1}} G \times \mc{D}_{x_i}^{\times} \xrightarrow{\tau_{x_i}} \mc{P}|_{\mc{D}_{x_i}^{\times}}.
\]
\item $\tau^{g}= (\tau_{x_i}^{g})_{i=1, \ldots, n}$ consists of the \textit{tautological} trivializations of $\mc{P}^g$ over each $\mc{D}_{x_i}$ coming from the definition of $\mc{P}^g$.
\end{itemize}
\end{defn}
At each of the points $x_i$, we now fix a pro-algebraic subgroup $\mbf{P}_i$ of $LG$ that is stable under the action of $\Aut_{\mc{O}}$; we additionally require that for some integer $m$, $\mbf{P}_i$ should contain the subgroup
\[
\mbf{I}(m)= \{g \in L^+G: g \equiv 1 \pmod {t^m}\}
\]
in finite co-dimension.
\begin{defn}\label{defnlevel}
Having fixed $S= \{x_1, \ldots, x_n\}$ and $\mbf{P}_1, \ldots, \mbf{P}_n$ as above, we define $\Bun_{G,S}(\mbf{P}_1, \ldots, \mbf{P}_n)$ to be the stack associated to the quotient prestack
\[
R \mapsto \Bun_{G, S, \infty}(R)/ \prod_{i=1}^n (\mbf{P}_i \rtimes \Aut_{\mc{O}})(R).
\]
When there is no risk of confusion, we omit the subscript $S$ from the notation and simply write $\Bun_G(\mbf{P}_1, \ldots, \mbf{P}_n)$.
\end{defn}
Note that since the action of $(LG \rtimes \Aut_{\mc{O}})^n$ does not necessarily preserve the isomorphism class of the $G$-torsor $\mc{P}$ on $X_R$, the moduli space $\Bun_{G, S}(\mbf{P}_1, \ldots, \mbf{P}_n)$ need not have a projection to $\Bun_G$. The action does not alter $\mc{P}|_{X_R- \cup \Gamma(x_i)}$, however, so an object of $\Bun_{G, S}(\mbf{P}_1, \ldots, \mbf{P}_n)(R)$ does yield a well-defined $G$-torsor on this complement. Also, the category $\Bun_{G, S}(\mbf{P}_1, \ldots, \mbf{P}_n)$ has a tautological object given by taking the image of an object of $\Bun_{G, S, \infty}^{\mr{pre}}(k)$ given by the trivial bundle with its tautological level structures and any fixed choice of local coordinates $\alpha_{x_i}$.\footnote{For any two such choices, the resulting objects of $\Bun_{G, S, \infty}(k)$ become uniquely isomorphic modulo the $\Aut_{\mc{O}}^n$-action.}
\begin{lemma}\label{algstack}
$\Bun_{G, S}(\mbf{P}_1, \ldots, \mbf{P}_n)$ is an algebraic stack locally of finite-type.
\end{lemma}
\proof
This follows exactly as in \cite[Corollary 4.2.6]{yun:exceptional}, by first deducing the result for 
\[
\Bun_{G, S}(\mbf{I}_1(m), \ldots, \mbf{I}_n(m))
\] 
from the (well-known) result for $\Bun_G$, and then deducing the case of $\Bun_{G, S}(\mbf{P}_1, \ldots, \mbf{P}_n)$ from that of $\Bun_{G, S}(\mbf{I}_1(m), \ldots, \mbf{I}_n(m))$.
\endproof
Just as in \cite[Lemma 4.2.5]{yun:exceptional}, we also have:
\begin{lemma}\label{omega}
For each $i=1, \ldots, n$, let 
\[
\Omega_{x_i}= N_{LG}(\mathbf{P}_i)/\mathbf{P}_i.
\]
Then there is a right-action of $\Omega_{x_i}$ on $\Bun_{G, S}(\mbf{P}_1, \ldots, \mbf{P}_n)$.
\end{lemma}
Finally, we can replace any $\mbf{P}_i$ by some finite cover, still acting on $\Bun_{G, S, \infty}$ on the right through $\mbf{P}_i$; Lemmas \ref{algstack} and \ref{omega} continue to hold.
\begin{rmk}
For the reader's convenience, we put this statement in its classical context: on automorphic forms over a function field $F$
\[
f \colon G(F) \bs G(\mathbb{A}_F)/\left( \prod_{x \neq x_i} G(\widehat{\mc{O}}_x) \times \prod_i \mbf{P}_i\right) \to \Qlb,
\]
we have the usual action by Hecke correspondences arising from decomposing the double coset $\mbf{P}_i w \mbf{P}_i$ into single cosets. But when $w$ normalizes $\mbf{P}_i$, the Hecke action comes from an actual automorphism (right-translation) of the moduli space $G(F) \bs G(\mathbb{A}_F)/\left( \prod_{x \neq x_i} G(\widehat{\mc{O}}_x) \times \prod_i \mbf{P}_i\right)$.
\end{rmk}
\section{Our setting}\label{setting}
Now we describe in detail the setting of this paper, taking \cite{yun:exceptional} as our starting-point. Let $G$ be a split (almost-)simple simply-connected group over $k$, satisfying the following two hypotheses:
\begin{itemize}
\item $G$ is oddly-laced;
\item $-1$ belongs to the Weyl group $W_G$ of $G$.
\end{itemize}
Explicitly, we take $G$ to be a split simple simply-connected group of type $A_1$, $D_{2n}$, $G_2$, $E_7$, or $E_8$ in the Dynkin classification. In fact, as we will see, the results of this paper are only non-trivial when the simply-connected and adjoint forms of $G$ differ: so for all practical purposes, we are working with types $A_1$, $D_{2n}$, and $E_7$.

Let $\tG$ be a split connected reductive group over $k$ with derived group equal to $G$, so that the quotient $\tG/G= S$ is a torus; call the quotient map $\nu \colon \tG \to S$. Fix a maximal torus $\tT$ of $\tG$ and a Borel $\tB$ containing $\tT$, likewise giving $T= \tT \cap G$, $B= \tB \cap G$, and determining based root data for $\tG$ and $G$, and an explicit Weyl group $W_G$ defined in terms of $T$. We denote by $\tZ$ and $Z_G$ the centers of $\tG$ and $G$, and we let $\tZ^0$ be the identity component of $\tZ$. Note that in all cases under consideration $Z_G= Z_G[2]$. The cases of particular interest for us--in which there is a non-trivial Kuga-Satake lifting problem--are those in which $Z_G \neq \{1\}$, namely types $A_1$, $D_{2n}$, and $E_7$. From now on we 
\begin{equation}\label{chark}
\text{assume the characteristic of $k$ is not 2.}
\end{equation}
In particular, $Z_G$ is a discrete group scheme over $k$, and the order of the kernel of the isogeny $\tZ \onto S$ is invertible in $k$. Our first task is to define the moduli spaces of $\tG$-bundles on $X= \mathbb{P}^1$ with level structure that will supply us with Hecke eigensheaves. We first recall the construction in \cite{yun:exceptional}. Yun works with the following conjugacy class of parahoric subgroups in $LG$ (see \cite[\S 2.2-2.3]{yun:exceptional}). In the apartment $\mc{A}(T)$ associated to $T$ of the building of $LG$, we can choose as origin the point corresponding to the subgroup $L^+G$, with the resulting identification $\mc{A}(T) \cong X_\bullet(T) \otimes \RR$. Then under this identification $\frac{1}{2} \rho^\vee$ lies in a unique facet, and we let $\mathbf{P}_{\frac{1}{2} \rho^\vee}$ be the parahoric subgroup associated to this facet. More precisely, Bruhat-Tits theory provides, for any facet $a$ in the building of $LG$, a smooth group scheme $\mc{P}_a$ over $k[[t]]$ with connected fibers whose generic fiber is $G \times_{\Spec k} \Spec k((t))$. We define $\mbf{P}_a$ to be the pro-algebraic subgroup of $LG$ representing the functor (of $k$-algebras)
\[
R \mapsto \mc{P}_a(R[[t]]).
\] 
We then apply this construction to the case where $a$ is the facet containing $\frac{1}{2} \rho^\vee$. Let $K$ denote the maximal reductive quotient of $\mbf{P}_{\frac{1}{2} \rho^\vee}$; since $G$ is simply-connected, $K$ is connected. Moreover, Yun shows (\cite[\S 2.5]{yun:exceptional}) that $K$ has a canonical connected double cover:
\begin{defn}
Let $\twoK$ denote the connected double-cover of $K$, so there is an exact sequence 
\[
1 \to \mu_2^{\mr{ker}} \to \twoK \to K \to 1.
\]
Note that our notation differs from that of \cite[\S 2.5]{yun:exceptional}, where this group is denoted $\widetilde{K}$; we reserve $\widetilde{(\ast)}$ for groups associated with the enlargement $\tG$ of $G$.
\end{defn}

We now define the particular moduli stacks of interest, beginning with the ones used in \cite{yun:exceptional}. Let $\mbf{P}_0 \subset L_0 G$ be the parahoric subgroup in the conjugacy class of $\mbf{P}_{\frac{1}{2} \rho^\vee}$ that contains the Iwahori $\mbf{I}_0 \subset L_0^+ G$, defined in terms of $B$. Moreover let  
\[
\twoP= \mbf{P}_0 \times_{K_0} \twoK_0,
\]
and let $\mbf{P}_0^+$ denote the pro-unipotent radical of $\mbf{P}_0$. Next let $\mbf{P}_{\infty}$ be the parahoric in the conjugacy class of $\mbf{P}_{\frac{1}{2} \rho^\vee}$ that contains the Iwahori $\mbf{I}_{\infty}^{\mr{op}} \subset L_{\infty}^+ G$ defined in terms of $B^{\mr{op}}$. Finally, let $\mbf{I}_1 \subset L_1^+ G$ denote the Iwahori subgroup defined again in terms of $B$. In the notation of \S \ref{levelstructure}, we now let $S= \{0, 1, \infty\} \subset \mathbb{P}^1(k)$; for later reference, we let $X^0$ be the variety $\mathbb{P}^1-S$ over $k$. The primary object of study in \cite{yun:exceptional} is the moduli space (see Definition \ref{defnlevel})
\[
\Bun= \Bun_G(\twoP, \mbf{I}_1, \mbf{P}_{\infty}).
\]
This sits in the following diagram:
\begin{equation}\label{basicdiagramG}
\xymatrix{
\Bun^+ \ar[d] \ar[r] & \Bun_G(\mbf{P}_0^+, \mbf{P}_{\infty}) \ar[d] \ar[dr] & \\
\Bun \ar[r] & \Bun_G(\twoP, \mbf{P}_{\infty}) \ar[r] & \Bun_G(\mbf{P}_0, \mbf{P}_{\infty}),
}
\end{equation}
where $\Bun^+= \Bun_G(\mbf{P}_0^+, \mbf{I}_1, \mbf{P}_{\infty})$. The vertical maps are $\twoK_0$-torsors, and the square is 2-cartesian.

Next we modify these constructions to define the corresponding moduli stacks of $\tG$-bundles on $X$. There are various ways of doing this; we take care to choose the new level structures so that the moduli spaces in the $G$ and $\tG$ cases are most easily compared.
\begin{defn}\label{tGlevels}
Let $\widetilde{\mbf{P}}_{\infty}$ be the sub-group scheme of $L_{\infty} \tG$ generated by $\mbf{P}_{\infty}$ and $L^+_{\infty} (\tZ^0)$. Let $\widetilde{\mbf{P}}_0(1)$ be the sub-group scheme of $L_0 \tG$ generated by $\mbf{P}_0$ and the pro-algebraic group $\tZ^0(1)$ defined as the kernel of reduction modulo $t$ (a local coordinate at zero),
\[
\tZ^0(1)= \ker( L^+_0(\tZ^0) \to \tZ^0).
\] 
\end{defn}
Note that $\widetilde{\mbf{P}}_0(1)$ is isomorphic to the direct product $\mbf{P}_0 \times \tZ^0(1)$: the restriction of $\nu$ to $\nu \colon \tZ^0 \to S$ can be identified with a product of maps $\mathbb{G}_m \to \mathbb{G}_m$, each given by multiplication by some $n \in \{\pm 1, \pm 2\}$, so (working one coordinate at a time) for any $x \in R[[t]]$, $1= \nu(1+tx)= 1+ntx+ (\text{higher order terms})$ forces $x=0$, since we have assumed (see (\ref{chark})) that 2 is invertible in $k$. Moreover, the pro-unipotent radical of $\widetilde{\mbf{P}}_0(1)$ is $\mbf{P}_0^+ \cdot \tZ^0(1)$, so the maximal reductive quotient of $\tPone$ is also $K_0$. In particular, we can form the analogous group ${}^{(2)}\tPone$ by pullback.

Finally, let $\widetilde{\mbf{I}}_1$ denote the Iwahori subgroup associated to $\widetilde{B}$ in $L_1^+ \tG$. With this notation in place, we introduce our main object of study:
\begin{defn}
Let $\tBun$ denote the algebraic stack $\Bun_{\tG}({}^{(2)}\tPone, \tI_1, \tPinf)$. 
\end{defn}
Similarly setting $\tBun^+= \Bun_{\tG}(\mbf{P}_0^+ \cdot \tZ^0(1), \tI_1, \tPinf)$, we then have the $\tG$-analogue of the basic diagram (\ref{basicdiagramG}):
\begin{equation}\label{basicdiagramtG}
\xymatrix{
\tBun^+ \ar[d] \ar[r] & \Bun_{\tG}(\mbf{P}_0^+ \cdot \tZ^0(1), \tPinf) \ar[d] \ar[dr] & \\
\tBun \ar[r] & \Bun_{\tG}({}^{(2)}\tPone, \tPinf) \ar[r] & \Bun_{\tG}(\tPone, \tPinf),
}
\end{equation}
in which the vertical maps are still $\twoK_0$-torsors, and the diagram is 2-cartesian.

We now must recall the Birkhoff decomposition and uniformization results for $G$- (or $\tG$-) bundles on $X= \mathbb{P}^1$. Consider the `trivial $G$-bundle on $\mathbb{A}^1$ with tautological $\mbf{P}_0$-level structure' $\mc{P}_{\mathbb{A}^1}^0$; to be precise, $\mc{P}^0_{\mathbb{A}^1}$ is defined as in \S \ref{levelstructure}, and is not literally a $G$-bundle on $\mathbb{A}^1$. Likewise let $\widetilde{\mc{P}}^0_{\mathbb{A}^1}$ be the trivial $\tG$-bundle on $\mathbb{A}^1$ with tautological $\tPone$ level structure at zero. Let $\Gamma_0$ and $\widetilde{\Gamma}_0$ denote the group ind-schemes of automorphisms of $\mc{P}_{\mathbb{A}^1}^0$ and $\widetilde{\mc{P}}^0_{\mathbb{A}^1}$, respectively. Also let $W^{\mr{aff}}$ denote the affine Weyl group $X_{\bullet}(T) \rtimes W_G$, and let $\widetilde{W}= X_{\bullet}(\tT) \rtimes W_G$ denote the Iwahori-Weyl group of $\tG$. The Weyl group of the reductive quotient $K_{\infty}$ of $\mbf{P}_{\infty}$ can be identified with a subgroup of $W^{\mr{aff}}$: take the subgroup generated by simple reflections that fix the alcove of $\mbf{P}_{\infty}$. The same holds for the reductive quotient $K_0$ of $\mbf{P}_0$ and its Weyl group, and in both cases, we write the resulting subgroup of $W^{\mr{aff}}$ as $W_K$.
\begin{lemma}\label{uniformization}
There are isomorphisms of stacks
\begin{align}
[\Gamma_0 \bs L_{\infty}G / \mbf{P}_{\infty}] &\xrightarrow{\sim} \Bun_G(\mbf{P}_0, \mbf{P}_{\infty}) \\
[\widetilde{\Gamma}_0 \bs L_{\infty} \tG / \tPinf] &\xrightarrow{\sim} \Bun_{\tG}(\tPone, \tPinf).
\end{align}
and Birkhoff decompositions
\begin{align}
L_{\infty}G (\bar{k})&= \coprod_{W_K \bs W^{\mr{aff}} / W_K} \Gamma_0(\bar{k}) w \mbf{P}_{\infty}(\bar{k})\\
L_{\infty} \tG(\bar{k}) &= \coprod_{W_K \bs \widetilde{W} / W_K} \widetilde{\Gamma}_0(\bar{k}) w \tPinf(\bar{k}).
\end{align}
\end{lemma} 
\proof
See \cite[\S 3.2.2]{yun:exceptional} and \cite[Proposition 1.1]{heinloth-ngo-yun:kloosterman}.
\endproof
It follows easily from diagram (\ref{basicdiagramtG}) that $\pi_0(\tBun)$ is naturally in bijection with $\pi_0(\Bun_{\tG}(\tPone, \tPinf))$; this is in turn in bijection (since $G$ is simply-connected) with 
\[
\pi_0(\Bun_{\tG}) \xrightarrow[\sim]{\nu} \pi_0(\Bun_S) \xleftarrow{\sim} \pi_0(L_{\infty} S/L_{\infty}^+ S) \xleftarrow{\sim} X_{\bullet}(S).
\]
We can describe the connected components of $\tBun$ in terms of this uniformization. First note that replacing $\mbf{P}_0$ with $\twoP$, and $\tPone$ with ${}^{(2)}\tPone$ we get obvious analogues of Lemma \ref{uniformization}. Then, for each $w \in W_K \bs \widetilde{W} / W_K$ we obtain an object $\widetilde{\mc{P}}_w$ of $\Bun_{\tG}({}^{(2)}\tPone, \tPinf)(k)$ by glueing $\widetilde{\mc{P}}^0_{\mathbb{A}^1}$ with $\Ad(w) \tPinf$; and we can make the corresponding construction of $\mc{P}_w \in \Bun_G(\twoP, \mbf{P}_{\infty})$ for $w \in W^{\mr{aff}}$. The stabilizers of $\mc{P}_w$ and $\widetilde{\mc{P}}_w$ are, respectively
\begin{align}
\Stab_w^G &= \left(\Gamma_0 \cap w \mbf{P}_{\infty} w^{-1}\right) \times_{K_0} \twoK_0 \\
\Stab_w^{\tG} &= \left( \widetilde{\Gamma}_0 \cap w \tPinf w^{-1} \right) \times_{K_0} \twoK_0.
\end{align}
In other words, $\Bun_G(\twoP, \mbf{P}_{\infty})$ has a stratification by sub-stacks $[\{\mc{P}_w\}/\Stab^G_w]$, and likewise $\Bun_{\tG}({}^{(2)}\tPone, \tPinf)$ has a stratification by sub-stacks $[\{\widetilde{\mc{P}}_w\}/\Stab^{\tG}_w]$. By taking the preimages in $\Bun$ and $\tBun$, we obtain stratifications by sub-stacks that we denote $\Bun_w$ (for $w \in W_K \bs W^{\mr{aff}} /W_K$) and $\tBun_w$ (for $w \in W_K \bs \widetilde{W} /W_K$), respectively. For $w= \lambda \rtimes w_G \in \widetilde{W}= X_{\bullet}(\tT) \rtimes W_G$, $\tBun_w$ lies in the component corresponding to $\nu \circ \lambda \in X_{\bullet}(S)$. In particular, we can identify the connected component of $\tBun$ containing the tautological object $\mc{P}_1$ as
\[
\tBun^0= \coprod_{\substack{w= \lambda \rtimes w_G \in W_K \bs \widetilde{W} / W_K \\
\lambda \in X_{\bullet}(T)}} \tBun_w= \coprod_{w \in W_K \bs W^{\mr{aff}}/ W_K} \tBun_w.
\]

Taking the associated $\tG$-bundle defines a map $\Bun \to \tBun$, and for $w \in W_K \bs W^{\mr{aff}} / W_K$ it respects the above stratifications, yielding a map $\Bun_w \to \tBun_w$. The crucial point is the following:
\begin{prop}\label{conndcomp}
The map $\Bun \to \tBun^0$ is an equivalence.
\end{prop}
\proof
We check this stratum-by-stratum. It suffices to show that for all $w \in X_{\bullet}(T) \rtimes W_G= W^{\mr{aff}} \subset \widetilde{W}$, $\Stab_w^G= \Stab_w^{\tG}$, i.e. that the natural map
\[
\Gamma_0 \cap w \mbf{P}_{\infty} w^{-1} \to \widetilde{\Gamma}_0 \cap w \tPinf w^{-1} 
\]
is an isomorphism. For a $k$-algebra $R$, an element of $\left( \widetilde{\Gamma}_0 \cap w \tPinf w^{-1} \right)(R)$ gives rise fppf-locally on $R$ to an equation of the form $p_0 z_0= w z_{\infty} p_{\infty} w^{-1}$ with $p_0 \in \mbf{P}_0(R)$, $z_0 \in \tZ^0(1)(R)$, $p_{\infty} \in \mbf{P}_{\infty}(R)$, and $z_{\infty} \in L^+_{\infty}(\tZ^0)$. Applying $\nu$, we find $\nu(z_0)= \nu(z_{\infty})$; but since $1+t R[[t]] \cap R[[t^{-1}]]^{\times}= \{1\}$, we see that $\nu(z_0)= \nu(z_{\infty})=1$. This forces (as in the argument following Definition \ref{tGlevels}, by our assumption on $\mr{char}(k)$) $z_0=1$, and $z_{\infty} \in \mbf{P}_{\infty}(R)$. We may as well then assume $z_{\infty}=1$ (incorporating $z_{\infty}$ into $p_{\infty}$), and so we actually have an equality $p_0= w p_{\infty} w^{-1}$ bearing witness to an element of $\left( \Gamma_0 \cap w \mbf{P}_{\infty} w^{-1}\right) (R)$. This implies 
\[
\Gamma_0 \cap w \mbf{P}_{\infty} w^{-1} \to \widetilde{\Gamma}_0 \cap w \tPinf w^{-1} 
\]
is an epimorphism, and as it is obviously injective, we are done.
\endproof
\section{The eigensheaves}
\subsection{Construction of the eigensheaves}\label{ramifiedHecke}
In this section, we combine the equivalence $\Bun \xrightarrow{\sim} \tBun^0$ of Proposition \ref{conndcomp} with the analysis of the sheaf theory of $\Bun$ carried out in \cite[Theorem 3.2]{yun:exceptional} to produce our desired Hecke eigensheaves on $\tBun$. The key simplification arises from applying Lemma \ref{omega} at the point $x=1$, where we have taken $\tI_1$ level structure. In this case we identify the group $\Omega_1$ with the stabilizer in $\widetilde{W}$ of the alcove corresponding to the standard Iwahori $\tI_1$, and $\nu \colon \Omega_1 \xrightarrow{\sim} X_{\bullet}(S)$ also identifies $\Omega_1$ with $\pi_0(\tBun)$. For $\gamma \in \Omega_1$, we denote by 
\[
\mathbb{T}_{\gamma} \colon \tBun \to \tBun
\]
the action given by Lemma \ref{omega}. Writing $\tBun^\gamma$ for the connected component corresponding to $\gamma$, we see that $\mathbb{T}_{\gamma}$ induces isomorphisms
\[
\mathbb{T}_{\gamma} \colon \tBun^0 \xrightarrow{\sim} \tBun^{\gamma}.
\]
In particular, all connected components of $\tBun$ are isomorphic (compare \cite[Corollary 1.2]{heinloth-ngo-yun:kloosterman}).\footnote{Note that this is a special, and highly simplifying, feature of our particular context: for contrast, observe that the degree zero and degree one connected components of $\Bun_{\mr{GL}_2}$ ($X= \mathbb{P}^1$ still) are not isomorphic, since no degree 1 vector bundle has $\mr{GL}_2$ as its automorphism group.} The idea is to take Yun's construction of a perverse Hecke eigensheaf on $\Bun \xrightarrow{\sim} \tBun^0$, and then use the (`ramified Hecke operators') $\mathbb{T}_{\gamma}$ to propagate the eigensheaf to the other connected components of $\tBun$. We begin by reviewing Yun's construction (\cite[\S 3]{yun:exceptional}). The tautological object in $\Bun_G(\mbf{P}_0, \mbf{P}_{\infty})$ (with automorphism group $K_0$) has preimage in $\Bun$ equivalent to a quotient $[\twoK_0 \bs fl_G]$, for a suitable action of $K_0$ on $fl_G$ (see \cite[\S 3.2.4]{yun:exceptional}). $K_0$ acts on $fl_G$ with finitely many orbits, so there is a unique open orbit $U \subset fl_G$, giving open embeddings
\[
[\twoK_0 \bs U] \subset [\twoK_0 \bs fl_G] \subset \Bun.
\]
As in \cite[\S 3.2.5]{yun:exceptional}, we fix a point $u_0 \in U(\Z[1/N])$ (for some $N$ sufficiently large, and for an integral model of $U$ arising from extending $K_0$ and $G$ to split reductive group schemes over some $\Z[1/M]$), and 
\begin{equation}\label{Upoint}
\text{denote by $u_0 \in U(k)$ the induced $k$-point for all $k$ of sufficiently large characteristic.}
\end{equation} 
This choice is in effect from now on. As an element of $U(k) \subset fl_G(k)$, $u_0$ corresponds to a Borel subgroup $B_0 \subset G$ over $k$, which is in general position with respect to $K_0$: 
\begin{defn}\label{thefinitegroup}
Let $\mr{A}$ denote the finite group scheme $B_0 \cap K_0$ over $k$. Let ${}^{(2)}\mathrm{A}$ denote the double-cover of $\mr{A}$ given by pullback along $\twoK_0 \onto K_0$.\footnote{Note that this is what Yun denotes $\widetilde{A}$.} Finally, let $Z({}^{(2)}\mathrm{A})$ denote the center of ${}^{(2)}\mathrm{A}$.
\end{defn}
Recall the following results (\cite[\S 2.6]{yun:exceptional}) on the structure and representation theory of the finite 2-group ${}^{(2)}\mathrm{A}(\bar{k})$. Recall that we have set
\begin{equation}
\Q'= \begin{cases}
\Q \quad \text{if $G$ is of type $D_{4n}$, $G_2$, or $E_8$;}\\
\Q(\imath) \quad \text{if $G$ is of type $A_1$, $D_{4n+2}$, or $E_7$,}
\end{cases}
\end{equation}
and have also set $\Ql'= \Ql(\imath)$. All sheaves considered will be $\Ql'$-sheaves. In parallel to this condition on the coefficients, we impose the following restriction on the field of definition $k$, in effect for the rest of this paper:
\begin{equation}\label{khassqrt}
\text{$\sqrt{-1} \in k$ for $G$ of type $A_1$, $D_{4m+2}$, or $E_7$.}
\end{equation}
\begin{lemma}\label{grouptheory}
Assume $k$ satisfies condition (\ref{khassqrt}), so that $\gal{k}$ acts trivially on $Z(\twoA)(\bar{k})$.
\begin{enumerate}
\item $\twoZ= Z(\twoA)$. 
\item Restriction to $\twoZ(\bar{k})$ gives a bijection between irreducible odd representations of $\twoA(\bar{k})$ and odd characters of $Z(\twoA(\bar{k}))$:
\begin{equation}\label{irrep}
\mr{Irr}_{\Qb}(\twoA(\bar{k}) )_{\mr{odd}} \xrightarrow{\sim} \Hom(Z(\twoA)(\bar{k}), \Qb^\times)_{\mr{odd}}= \Hom(Z(\twoA)(\bar{k}), \Q'^\times)_{\mr{odd}}.
\end{equation}
\item If $k$ is a finite field, local field, or number field, then for each odd $\chi \colon Z(\twoA)(\bar{k}) \to \Q'^\times$, the corresponding irreducible representation $V_{\chi}$ of $\twoA(\bar{k})$ descends to an irreducible representation of $\twoA(\bar{k}) \rtimes \gal{k}$, whose coefficients can be taken to be $\Q(\imath)$.
\end{enumerate}
\end{lemma}
\proof
The first claim is \cite[Lemma 2.6(2)]{yun:exceptional}. The second claim is elementary: the inverse of the isomorphism (\ref{irrep}) is given by inducing the central character, up to some multiplicity. The third claim is a variant of \cite[Lemma 2.7]{yun:exceptional}, whose proof is not complete.\footnote{Namely, that argument uses the incorrect assertion that $H^2(\overline{\Gamma}, \Qb^\times)=0$ for $\overline{\Gamma}$ a finite group isomorphic to a direct sum of $\Z/2\Z$'s.} The obstruction to descending $V_{\chi}$ to a representation of $\twoA(\bar{k}) \rtimes \gal{k}$ is a class in $H^2(\gal{k}, \Qb^\times)$. This Galois cohomology group vanishes for the claimed $k$: this is elementary for $k$ finite, and for local and especially number fields it is a beautiful theorem of Tate (\cite[Theorem 4]{serre:DSsurvey}). The argument showing the descended $V_{\chi}$ can be defined with $\Q(\imath)$ coefficients is as in \cite[Lemma 2.7]{yun:exceptional}.
\endproof
For clarity, we collect in one place the various conditions in effect on the field of definition $k$:
\begin{defn}\label{Vchi}
Consider any odd central character $\chi \colon Z(\twoA)(\bar{k}) \to \Q'^\times$, with associated irreducible representation $V_{\chi}$ of $\twoA(\bar{k})$. Let $k$ be any field satisfying conditions (\ref{khassqrt}) and (\ref{Upoint}), and \textit{moreover} for which $V_{\chi}$ satisfies the conclusion of Lemma \ref{grouptheory}(3). Then from now on let $V_{\chi}$ denote a fixed choice of descent to an irreducible representation of $\twoA(\bar{k}) \rtimes \gal{k}$, with $\Q(\imath)$ coefficients.
\end{defn}

We now recall the crucial result analyzing the sheaf theory of $\Bun$, or, in our case, $\tBun^0$. Throughout, for an algebraic stack $\mf{X}$ over a field $k$, we will write $D^b(\mf{X})$ for the derived category of bounded complexes of $\Ql'$-sheaves with constructible cohomology, as in \cite{laszlo-olsson:sixops2} (if we need to specify another field of coefficients, $\Ql$ for instance, we will write $D^b(\mf{X}, \Ql)$). Recall (\cite[\S 3.3.1]{yun:exceptional}) the sub-category 
\[
D^b(\Bun)_{\mr{odd}} \subset D^b(\Bun)
\]
of odd sheaves, on which $\mu_2^{\mr{ker}}= \ker(\twoK_0 \to K_0)$ acts by the sign character. We can similarly define $D^b(\tBun)_{\mr{odd}}$, since $\mu_2^{\mr{ker}}$ is also contained in the automorphism group of every object of $\tBun$. For future reference, let us also note a refinement of this observation: the automorphism group of every object of $\Bun_{\tG}(\tPone, \tI_1, \tPinf)$ contains the center $Z_G$ of $G$, and likewise the automorphism group of every object of $\tBun$ contains the double cover (pullback under $\twoK_0 \to K_0$) ${}^{(2)}Z_G$ of $Z_G$. We can therefore decompose $D^b(\tBun)$ into a direct sum of categories $D^b(\tBun)_{\psi}$, indexed over characters $\psi \colon {}^{(2)}Z_G \to \Qlb^\times$. We of course will be interested in the corresponding decomposition of $D^b(\tBun)_{\mr{odd}}$ into a direct sum over the odd characters $\psi$. 

Now we recall the main result analyzing odd sheaves on $\Bun$. Let $j \colon [\twoK_0 \bs U] \into \Bun$ denote the open inclusion.
\begin{thm}[Theorem 3.2 of \cite{yun:exceptional}]
Assume $G$ is the split simple simply-connected group of type $A_1$, $D_{2n}$, $E_7$, $E_8$, or $G_2$. Then the restriction
\[
j^* \colon D^b(\Bun)_{\mr{odd}} \to D^b([\twoK_0 \bs U])_{\mr{odd}}
\]
is an equivalence of categories with quasi-inverse given by $j_!= j_*$.
\end{thm}
The analysis of connected components of $\tBun$ then implies:
\begin{cor}\label{sheaftheory}
For all $\gamma \in \Omega_1$, consider the composite
\[
j_{\gamma}= \mathbb{T}_{\gamma} \circ j \colon [\twoK_0 \bs U] \into \tBun^{\gamma}.
\]
Then the restriction
\[
j_{\gamma}^* \colon D^b(\tBun^{\gamma})_{\mr{odd}} \to D^b([\twoK_0 \bs U])_{\mr{odd}}
\]
is an equivalence with inverse $j_{\gamma, !}= j_{\gamma, *}$.
\end{cor}
Assume $k$ is as in Definition \ref{Vchi}. We can now define the hoped-for eigensheaves on $\tBun$ over $k$, starting from Yun's construction on $\Bun$. Fix an odd character (recall equation (\ref{irrep}))
\begin{equation}\label{chi}
\chi \colon Z({}^{(2)}\mathrm{A})(\bar{k}) \to \Q'^\times, 
\end{equation}
to which we have associated (Lemma \ref{grouptheory} and Definition \ref{Vchi}) an irreducible representation $V_{\chi}$ of ${}^{(2)}\mathrm{A}(\bar{k}) \rtimes \gal{k}$ having $\chi$ as central character. By \cite[Lemma 3.3]{yun:exceptional}, $V_{\chi} \otimes_{\Q(\imath)} \Ql'$ is isomorphic to the pullback under $u_0$ of a geometrically irreducible local system
\[
\mc{F}_{\chi} \in \Loc_{\twoK_0}(U, \Ql')_{\mr{odd}}, 
\] 
which we view as an object of $D^b([\twoK_0 \bs U])_{\mr{odd}}$. Yun's eigensheaf is then (\cite[Theorem 4.2]{yun:exceptional})
\[
j_!(\mc{F}_{\chi})= j_*(\mc{F}_{\chi}) \in D^b(\Bun)_{\mr{odd}}.
\]
\begin{defn}
Assume $k$ is as in Definition \ref{Vchi}. Let $\chi \colon Z({}^{(2)}\mathrm{A}(\bar{k})) \to \Qlb^\times$ be any odd character. We let $\mc{A}_{\chi} \in D^b(\tBun)_{\mr{odd}}$ be the perverse sheaf on $\tBun$ whose restriction $\mc{A}_{\chi}^{\gamma}$, for all $\gamma \in \Omega_1$, to $\tBun^{\gamma}$ is given by
\[
\mc{A}_{\chi}^{\gamma}= \mc{A}_{\chi}|_{\tBun^{\gamma}}= j_{\gamma, !} \mc{F}_{\chi}= j_{\gamma, *} \mc{F}_{\chi}.
\]
\end{defn}
That is, we make the only definition compatible with the requirement that $\mc{A}_{\chi}^0$ be Yun's eigensheaf, and that $\mc{A}_{\chi}$ be eigen for the ramified Hecke operators $\mathbb{T}_{\gamma}$ at $1 \in X$.
\subsection{Geometric Satake equivalence}\label{satake}
We recall a convenient form of the geometric Satake equivalence. See \cite{mirkovic-vilonen:geometricsatake} and \cite[\S 4.1]{yun:exceptional} for more background. Let $\mc{G}$ be any split connected reductive group over $k$ ($\mc{G}$ will of course eventually be either $G$ or $\tG$). Let $\Gr_{\mc{G}}= L\mc{G}/L^+\mc{G}$ as usual denote the affine Grassmannian of $\mc{G}$. The main result of \cite{mirkovic-vilonen:geometricsatake} describes the category $\Sat^{\mr{geom}}_{\mc{G}}$ of $(L^+\mc{G})_{\bar{k}}$-equivariant perverse sheaves on $\Gr_{\mc{G}, \bar{k}}$: $\Sat^{\mr{geom}}_{\mc{G}}$ admits a convolution product making it a neutral Tannakian category over $\Ql$ with fiber functor
\begin{align}
H^* \colon \Sat^{\mr{geom}}_{\mc{G}} & \to \Vect_{\Ql} \\
\mc{K} &\mapsto H^*(\Gr_{\mc{G}, \bar{k}}, \mc{K})
\end{align}
This fiber functor induces an equivalence 
\[
\Sat^{\mr{geom}}_{\mc{G}} \xrightarrow{\sim} \Rep(\mc{G}^\vee)
\]
where we write $\mc{G}^\vee$ for the (split form over $\Ql$ of the) dual group of $\mc{G}$. We need a version of $\Sat^{\mr{geom}}_{\mc{G}}$ over $k$ rather than $\bar{k}$. It is natural for us to deviate from \cite[\S 4.1]{yun:exceptional} and instead follow the suggestion of \cite[Remark 2.9]{heinloth-ngo-yun:kloosterman} and \cite[\S 2]{frenkel-gross:rigid}. Recall that the simple objects of $\Sat^{\mr{geom}}_{\mc{G}}$ are given by the intersection cohomology sheaves of the affine Schubert varieties $\Gr_{\mc{G}, \leq \lambda}$. For all dominant $\lambda \in X_{\bullet}(T)$, we write
\[
j_{\lambda} \colon \Gr_{\mc{G}, \lambda} \into \Gr_{\mc{G}}
\]
for the inclusion of the $L^+ \mc{G}$-orbit containing $t^{\lambda}$. Then by definition the intersection cohomology sheaf of the closure $\Gr_{\mc{G}, \leq \lambda}$ of $\Gr_{\mc{G}, \lambda}$ is
\[
\IC_{\lambda}= j_{\lambda, !*} \Ql[\langle 2 \rho, \lambda \rangle],
\]
the shift reflecting that the dimension of $\Gr_{\mc{G}, \lambda}$ is $\langle 2 \rho, \lambda \rangle$. We will define $\Sat_{\mc{G}}$ to be the full subcategory of perverse sheaves on $\Gr_{\mc{G}}$ consisting of finite direct sums of arbitrary Tate twists $\IC_{\lambda}(m)$, for all $\lambda \in X_{\bullet}(T)^+$ and $m \in \Z$. Note that, in contrast to \cite[\S 4.1]{yun:exceptional}, we do not normalize the weights of the $\IC_{\lambda}$'s to be zero: this bookkeeping device frees us from having to choose a square root of the cyclotomic character;\footnote{Which of course cannot be done over $k= \Q$, although it is possible over many quadratic extensions of $\Q$.} and it ensures that the local systems we eventually construct will specialize (at points of $X^0(K)$, for $K/\Ql$ finite) to de Rham Galois representations. Adapting the argument of \cite[\S 4.1]{yun:exceptional} to our normalization, a result of Arkhipov-Bezrukavnikov (\cite[\S 3]{arkhipov-bezrukavnikov:pervonflag}) implies that $\Sat_{\mc{G}}$ is closed under convolution: to be precise, we have
\[
\IC_{\lambda} \ast \IC_{\mu} \cong \oplus_{\nu} \IC_{\nu}(\langle \nu- \lambda- \mu, \rho \rangle).
\]
Note that $\langle \nu- \lambda-\mu, \rho \rangle$ is an integer, since only $\nu$ for which there is an inclusion of highest weight representations $V_{\nu} \into V_{\mu} \otimes V_{\lambda}$, and in particular for which $\lambda+\mu-\nu$ lies in the root lattice, will appear on the right-hand side. We would like to combine the tensor functor
\begin{equation}
H^*_{\bar{k}} \colon \Sat_{\mc{G}} \to \Sat^{\mr{geom}}_{\mc{G}} \xrightarrow{H^*} \Rep(\mc{G}^\vee)
\end{equation}
with a mechanism for keeping track of the weight/Tate twist. One way to do this is to replace $\Sat_{\mc{G}}$ by a skeleton whose objects are precisely the direct sums of the various $\IC_{\lambda}(n)$ (the skeleton can be equipped with a suitable tensor structure making these equivalent as tensor categories), and then to define a fully faithful tensor functor
\[
H^*_w \colon \Sat_{\mc{G}} \to \Rep(\mc{G}^\vee \times \mathbb{G}_m)
\]
by additively extending the assignment on simple objects
\[
\IC_{\lambda}(n) \mapsto H^*_{\bar{k}}(\IC_{\lambda}(n)) \boxtimes \left( z \mapsto z^{\langle 2\rho, \lambda \rangle- 2n} \right).
\]
Composing with the canonical fiber functor $\omega$ of $\Rep(\mc{G}^\vee \times \mathbb{G}_m)$, this yields a surjective homomorphism $\mc{G}^\vee \times \mathbb{G}_m \to \Aut^{\otimes}(\omega \circ H^*_w)$ whose kernel 
\[
\{(g, z) \in \mc{G}^\vee \times \mathbb{G}_m: \text{for all dominant $\lambda \in X_{\bullet}(T)$ and all $n \in \Z$, $g$ acts on $V_{\lambda}$ by $z^{2n- \langle 2\rho, \lambda \rangle}$} \}
\]
is clearly equal to the subgroup $\langle (2\rho(-1), -1) \rangle \subset \mc{G}^\vee \times \mathbb{G}_m$. That is, we have a tensor-equivalence $\Sat_{\mc{G}} \xrightarrow{\sim} \Rep(\mc{G}_1^\vee)$, where we define (following \cite{frenkel-gross:rigid})
\begin{equation}\label{G1}
\mc{G}_1^\vee = (\mc{G}^\vee \times \mathbb{G}_m)/ \langle (2\rho(-1), -1)\rangle.
\end{equation}
Note that If $\mc{G}$ is simply-connected, then $\mc{G}_1^\vee$ is isomorphic to $G^\vee \times \mathbb{G}_m$, since $2\rho(-1)=1$.
\subsection{Geometric Hecke operators}
We briefly recall the definition of geometric Hecke operators in our context, as well as the notion of Hecke eigensheaf. Recall that the Hecke stack $\tHk$ associated to $\tBun$ is the category of tuples $(R, x, \mc{P}, \mc{P}', \iota)$ where:
\begin{itemize}
\item $R$ is a $k$-algebra;
\item $x \in X^0(R)$;
\item $\mc{P}$ and $\mc{P}'$ are objects of $\tBun(R)$;
\item $\iota$ is an isomorphism of $\mc{P}$ and $\mc{P}'$ away from the graph of $x$.
\end{itemize}
Projecting such data to $(R, x, \mc{P})$ (the map $\lefth$) or $(R, x, \mc{P}')$ (the map $\righth$) gives a correspondence diagram
\[
\xymatrix{
& \tHk \ar[dl]_{\overleftarrow{h}} \ar[dr]^{\overrightarrow{h}} & \\
\tBun \times X^0 & & \tBun \times X^0 \\
}
\]
As explained in \cite[\S 4.1.3]{yun:exceptional} (using the fact that--see \cite[Remark 4.1]{heinloth-ngo-yun:kloosterman}--$\righth$ and $\lefth$ are locally trivial fibrations in the smooth topology, with fibers isomorphic to $\Gr_{\tG}$), or slightly differently in \cite[4.3.1]{yun:iccmsurvey}, for each $\mc{K} \in \Sat_{\tG}$ there is an object $\mc{K}_{\tHk} \in D^b(\tHk, \Ql)$ whose restriction to each geometric fiber of $\righth$ is isomorphic to $\mc{K}$. As usual, the (universal) geometric Hecke operator is the functor
\begin{align}
\mathbb{T} \colon \Sat_{\tG} \times D^b(\tBun \times X^0) &\to D^b(\tBun \times X^0) \\
(\mc{K}, \mc{F}) &\mapsto \righth_! \left( \lefth^* (\mc{F}) \otimes_{\Ql} \mc{K}_{\tHk} \right)
\end{align}
The induced functor
\[
\Sat_{\tG} \to \End(D^b(\tBun \times X^0))
\]
is monoidal. When the input from $D^b(\tBun \times X^0)$ is of the form $\mc{F} \boxtimes \Qlb$ for some $\mc{F} \in D^b(\tBun)$, we write
\[
\mathbb{T}_{\mc{K}}(\mc{F})= \mathbb{T}(\mc{K}, \mc{F} \boxtimes \Qlb)
\]
Finally, recall the definition of a Hecke eigensheaf:
\begin{defn}\label{eigensheafdefn}
Let $\mc{F}$ be an object of $D^b(\tBun)$. We say that $\mc{F}$ is a Hecke eigensheaf if there exists
\begin{itemize}
\item a tensor functor $\tE \colon \Sat_{\tG} \to \Loc(X^0)$;
\item a system of isomorphisms $\epsilon_{\mc{K}}$ for all $\mc{K} \in \Sat_{\tG}$
\[
\epsilon_{\mc{K}} \colon \mathbb{T}_{\mc{K}}(\mc{F}) \xrightarrow{\sim} \mc{F} \boxtimes \tE(\mc{K})
\]
satisfying compatibility conditions that will not concern us (see \cite[following Proposition 2.8]{gaitsgory:dejongconjecture}).
\end{itemize}
In this case we call $\tE$ the eigen-local system of $\mc{F}$.
\end{defn}
\subsection{Proof of the eigensheaf property}
Recall that we have fixed a point $u_0 \colon \Spec k \to U$. We also write $u_0$ for the induced maps $\Spec k \to [\twoK_0 \bs U] \subset \tBun^0 \subset \tBun$. For all $\gamma \in \Omega_1$, we can compose with $\mathbb{T}_{\gamma}$ to obtain
\[
u_{\gamma} \colon \Spec k \to \tBun^{\gamma}.
\]
From Corollary \ref{sheaftheory}, we obtain equivalences
\[
(u_{\gamma} \times \id)^* \colon D^b(\tBun^{\gamma} \times X^0)_{\mr{odd}} \xrightarrow{\sim} D_{{}^{(2)}\mathrm{A}}(X^0)_{\mr{odd}},
\]
where ${}^{(2)}\mathrm{A}$ acts trivially on $X^0$. The strategy for proving $\mc{A}_{\chi}$ is an eigensheaf ($\chi$ as in equation \ref{chi}) is to show that for all $\gamma \in \Omega_1$ and all $\mc{K} \in \Sat_{\tG}$, $(u_{\gamma} \times \id)^* \mathbb{T}_{\mc{K}}(\mc{A}_{\chi})$ is concentrated in a single perverse degree. Such sheaves $\mc{A}_{\chi}$ can then be explicitly described via Corollary \ref{sheaftheory} and an analogue of \cite[Lemma 3.4]{yun:exceptional}. In preparation for this computation, note that the $\mathbb{T}_{\gamma, !}$ and $\mathbb{T}_{\gamma}^*$ commute with the $\mathbb{T}_{\mc{K}}$: informally this is the statement that `Hecke operators at different places commute'; more formally, the stack $\tHk$ carries an $\Omega_1$-action compatible with its two projections $\lefth$ and $\righth$. Furthermore, the spread-out sheaves $\mc{K}_{\tHk}$ (for all $\mc{K} \in \Sat_{\tG}$) are $\Omega_1$-equivariant, so we find
\begin{align}
(u_{\gamma} \times \id)^* \mathbb{T}_{\mc{K}}(\mc{A}_{\chi}) &\cong (u_0 \times \id)^*(\mathbb{T}_{\gamma} \times \id)^* \mathbb{T}_{\mc{K}}(\mc{A}_{\chi}) \\
&\cong (u_0 \times \id)^* \mathbb{T}_{\mc{K}}(\mathbb{T}_{\gamma}^* \mc{A}_{\chi}) \cong (u_0 \times \id)^* \mathbb{T}_{\mc{K}}(\mc{A}_{\chi}).
\end{align}
Now consider the following diagram, where declaring the squares cartesian defines the new objects $\tGR$ and $\tGR^U_{\gamma}$:
\begin{equation}\label{maindiagram}
\xymatrix{
& \tGR^U_{\gamma} \ar[ld]_{\omega^U_{\gamma}} \ar[r]^{j_{\gamma}} & \tGR \ar[ld]^{\omega} \ar[r] \ar[d]^{\pi} & \tHk \ar[d]^{\righth} \\
[\twoK_0 \bs U] \ar[r]^-{j_{\gamma}} & \tBun & X^0 \ar[r]^-{u_0 \times \id} & \tBun \times X^0.
}
\end{equation}
Here $\omega$ is the remaining projection corresponding to $\lefth$ on $\tHk$. Note that $\tGR$ is the analogue of the Beilinson-Drinfeld Grassmannian in this context.\footnote{Note that we continue to adhere to the notational pattern of using $\widetilde{(\ast)}$ to denote the $\tG$-version of an object that could similarly be defined for $G$. Our notation is as a result not always consistent with that of \cite{yun:exceptional}: for instance, $\tGR^U$ denotes there (the version for $G$ of) what we will call $\mf{G}^U$ below (see diagram (\ref{maindiagrambis})).} Let us also denote by
\[
\pi^U_{\gamma} \colon \tGR^U_{\gamma} \to X^0
\]
the composite $\pi \circ j_{\gamma}$. Repeated application of proper base-change yields
\begin{equation}\label{basechange}
(u_0 \times \id)^* \mathbb{T}_{\mc{K}}(j_{\gamma, !} \mc{F}_{\chi})= (u_0 \times \id)^* \righth_! \left( \lefth^*(j_{\gamma, !} \mc{F}_\chi) \otimes \mc{K}_{\tHk} \right) \cong \pi^U_{\gamma, !} \left( \omega_{\gamma}^{U, *}(\mc{F}_{\chi}) \otimes \mc{K}_{\tGR} \right),
\end{equation}
where $\mc{K}_{\tGR}$ denote the pull-back of $\mc{K}_{\tHk}$ to $\tGR$. $\mc{K}_{\tGR}[1]$ is perverse (recall the fibers of $\mc{K}_{\tGR}$ at $x \in X^0$ are copies of $\mc{K}$), and $\mc{F}_{\chi}$ is a local system (in cohomological degree zero), so $\omega_{\gamma}^{U, *}(\mc{F}_{\chi}) \otimes \mc{K}_{\tGR}[1]$ is perverse.    Our immediate aim is to show that each $(u_0 \times \id)^* \mathbb{T}_{\mc{K}}(j_{\gamma, !} \mc{F}_{\chi})[1]$ is a perverse sheaf on $X^0$. Any object $\mc{K}$ of $\Sat_{\tG}$ is a direct sum of simple objects, so we may assume $\mc{K}$ is simple and therefore supported on some  $\Gr_{\tG, \leq \lambda}$, $\lambda \in X_{\bullet}(\tT)$. The corresponding $\mc{K}_{\tGR}$ is then supported on a corresponding sub-stack $\tGR_{\leq \lambda}$, which pulls back in diagram (\ref{maindiagram}) to a substack $\tGR^U_{\gamma, \leq \lambda}$ of $\tGR^U_{\gamma}$.

We now come to the crucial geometric lemma. We note that Yun has found (see \cite[Lemma 4.4.7]{yun:CDM}) an argument that applies much more generally; the following, an elaboration of \cite[Lemma 4.8]{yun:exceptional} will suffice for us.
\begin{lemma}\label{affine}
For all $\gamma \in \Omega_1$, $\pi^U_{\gamma} \colon \tGR^U_{\gamma, \leq \lambda} \to X^0$ is affine.
\end{lemma}
\proof
Since $[\twoK_0 \bs U] \subset [\twoK_0 \bs fl_G]$ is affine, we may replace $\tGR^U_{\gamma, \leq \lambda}$ with the preimage of 
\[
[\twoK_0 \bs fl_G] \xrightarrow{j_{\gamma}} \tBun^{\gamma}.
\]
Let us call this $\tGR^{fl}_{\gamma, \leq \lambda}$. By construction as the preimage of $\mathbb{B}(K_0) \subset \Bun_{\tG}(\tPone, \tPinf)$ (under the morphism (\ref{projection}) below), and using \cite[Lemma 3.1]{yun:exceptional}, $\tGR^{fl}_{\gamma}$ (respectively $\tGR^{fl}_{\gamma, \leq \lambda}$) is the non-vanishing locus of a non-zero section $s$ of a line bundle $\mc{L}$ on $\tGR_{\gamma}$ (respectively $\tGR_{\gamma, \leq \lambda}$). It suffices to show the line bundle in question is ample. By \cite[Proposition 1.7.8]{lazarsfeld:positivity1}, this can be checked on geometric fibers, since the morphism $\tGR_{\gamma, \leq \lambda} \to X^0$ is proper. Thus let $x \colon \Spec K \to X^0$ be a geometric point of $X^0$, and consider the section $x^*s$ of $x^*\mc{L}$. The fiber $\tGR_{\gamma, \leq \lambda,x}$ is isomorphic to the $\gamma$-component, truncated by $\lambda$ of the affine Grassmannian $\Gr_{\tG}$; we denote this by $\Gr^{\gamma}_{\tG, \leq \lambda}$. We claim that $x^* \mc{L}$ is ample on $\Gr^{\gamma}_{\tG}$, so in particular its restriction to the closed sub-scheme $\Gr^{\gamma}_{\tG, \leq \lambda}$ is ample. This claim results from the following two assertions:
\begin{itemize}
\item $\Pic(\Gr_{\tG}^{\gamma}) \cong \Z$;
\item $x^*s$ is a non-zero global section of $x^* \mc{L}$ (which by the previous item must then be ample).
\end{itemize}
The first item follows from \cite[Corollary 12]{faltings:algloop}. To be precise, that result shows that $\Pic(\Gr_G) \cong \Z$ (for $G$ our simply-connected group), but the same then follows for each connected component of $\Gr_{\tG}$.\footnote{To be absolutely precise: consider along with the affine Grassmannian the affine flag variety $\Fl_{\tG}= L\tG/\widetilde{\mbf{I}}$, where $\widetilde{\mbf{I}}$ denotes the Iwahori. The connected components $\Fl_{\tG}^0$ and $\Gr_{\tG}^0$ are, up to taking reduced sub-schemes, isomorphic to their semi-simple counterparts $\Fl_G$ and $\Gr_G$ (see, eg, \cite[Proposition 6.6]{pappas-rapoport:twistedloop}). As in \S \ref{ramifiedHecke}, the different component of $\Fl_{\tG}$ are isomorphic via ramified Hecke operators $\Fl_{\tG}^0 \xrightarrow{\mathbb{T}_{\gamma}} \Fl_{\tG}^{\gamma}$. $\Pic(\Gr_{\tG}^0)$ is isomorphic to the subgroup of $\Pic(\Fl_{\tG}^0)$ corresponding to the unique minimal parahoric $\mbf{P}$ properly containing $\widetilde{\mbf{I}}$ but not contained in $L^+ \tG$ (see the proof of \cite[Corollary 12]{faltings:algloop}); by the same argument, $\Pic(\Gr_{\tG}^{\gamma})$ can be described inside of $\Pic(\Fl_{\tG}^{\gamma})$ as the subspace spanned by the natural $\mc{O}(1)$ on $\mathbb{T}_{\gamma}(\mbf{P})/\widetilde{\mbf{I}}$.} 
For the second item, recall that the pair $(\mc{L}, s)$ is the pull-back along the composite 
\begin{equation}\label{projection}
\tGR_{\gamma} \to \tBun^{\gamma} \xrightarrow[\sim]{\mathbb{T}_{\gamma}^{-1}} \tBun^0 \to \Bun_{\tG}(\tPone, \tPinf)^0 \xleftarrow{\sim}\Bun_{G}(\mbf{P}_0, \mbf{P}_{\infty}),
\end{equation}
where the original section is non-vanishing on the locus $\mathbb{B}K_0 \subset \Bun_G(\mbf{P}_0, \mbf{P}_{\infty})$ corresponding to the tautological object. It suffices then to show that the geometric fibers of $\tGR_{\gamma}$ over $\mathbb{B}K_0 \times X^0$ are non-empty. To see this, note that $\tHk \to \tBun \times X^0$ has geometric fibers isomorphic to $\Gr_{\tG}$. Choosing an element $\mc{P}$ of the fiber over $(\mc{P}_{u_0}, x)$ that lies in the $\gamma$ component of $\Gr_{\tG}$, we are done: the isomorphism $\iota \colon \mc{P}|_{X-\{x\}} \xrightarrow{\sim} \mc{P}_{u_0}|_{X-\{x\}}$ automatically implies that $\mc{P}$ projects to an object isomorphic to the tautological object of $\Bun_G(\mbf{P}_0, \mbf{P}_{\infty})$.  
\endproof
With Lemma \ref{affine} in hand, we can prove the main result of this section:
\begin{thm}\label{eigensheaf}
For all odd characters $\chi \colon Z({}^{(2)}\mathrm{A}) \to \Qlb^\times$, $\mc{A}_{\chi}$ is a Hecke eigensheaf.
\end{thm}
\proof
 Since $\omega_{\gamma}^{U, *}(\mc{F}_{\chi}) \otimes \mc{K}_{\tGR}[1]$ is perverse, and $\pi^U_{\gamma}$ is affine,
\[
(u_0 \times \id)^* \mathbb{T}_{\mc{K}}(j_{\gamma, !} \mc{F}_{\chi}) \cong \pi^U_{\gamma, !}\left(\omega_{\gamma}^{U, *}(\mc{F}_{\chi}) \otimes \mc{K}_{\tGR} \right) \in {}^p D^{\geq 1}(X^0).
\]
But by Corollary \ref{sheaftheory}, this is also
\begin{equation}\label{semiperverse}
(u_0 \times \id)^* \mathbb{T}_{\mc{K}}(j_{\gamma, *} \mc{F}_{\chi}) \cong \pi_! \left(\omega^* j_{\gamma, *}(\mc{F}_{\chi}) \otimes \mc{K}_{\tGR} \right).
\end{equation}
There is a natural isomorphism $\omega^* \circ j_{\gamma, *} \xrightarrow{\sim} j_{\gamma, *} \circ \omega^{U, *}_{\gamma}$: as in the proof of \cite[Proposition 4.7]{yun:exceptional}, this follows from the fact that $\lefth$ is a locally trivial fibration in the smooth topology. Thus, identifying $\pi_!= \pi_*$ on the support of $\mc{K}_{\tGR}$ ($\pi \colon \tGR_{\leq \lambda} \to X^0$ is proper), and using the projection formula and the Leray spectral sequence, we can carry on the identification \ref{semiperverse} as
\begin{equation}\label{dualsemiperverse}
(u_0 \times \id)^* \mathbb{T}_{\mc{K}}(j_{\gamma, !} \mc{F}_{\chi}) \cong \pi_* \left( (j_{\gamma, *} \omega^{U, *}_{\gamma} \mc{F}_{\chi}) \otimes \mc{K}_{\tGR} \right) \cong \pi^U_{\gamma, *} \left( \omega^{U, *}_{\gamma} \mc{F}_{\chi} \otimes \mc{K}_{\tGR} \right).
\end{equation}
(This is just the obvious variant of \cite[4.19]{yun:exceptional}.) Since $\pi^U_{\gamma}$ is affine, we can dually conclude that
\[
(u_0 \times \id)^* \mathbb{T}_{\mc{K}}(j_{\gamma, !} \mc{F}_{\chi}) \in {}^p D^{\leq 1}(X^0),
\]
hence that $(u_0 \times \id)^* \mathbb{T}_{\mc{K}}(j_{\gamma, !} \mc{F}_{\chi})[1]$ is perverse. Consequently, $(u_0 \times \id)^* \mathbb{T}_{\mc{K}}(\mc{A}_{\chi})[1]$ is perverse.

Now, for each component $\tBun^{\gamma}$ of $\tBun$, we apply \cite[Lemma 3.4]{yun:exceptional} to $(u_{\gamma} \times \id)^* \mathbb{T}_{\mc{K}}(\mc{A}_{\chi})$ to conclude
\begin{align}
\mathbb{T}_{\mc{K}}(\mc{A}_{\chi})|_{\tBun^{\gamma} \times X^0} &\cong \bigoplus_{\substack{\psi \colon Z({}^{(2)}\mathrm{A}) \to \Qlb^\times:\\ \text{$\psi$ is odd}}} (j_{\gamma, !} \mc{F}_{\psi}) \boxtimes \left(V_{\psi}^* \otimes (u_{\gamma} \times \id)^* \mathbb{T}_{\mc{K}}(\mc{A}_{\chi})_{\psi}\right)^{{}^{(2)}\mathrm{A}} \\
&= (j_{\gamma, !} \mc{F}_{\chi}) \boxtimes \left(V_{\chi}^* \otimes (u_{\gamma} \times \id)^* \mathbb{T}_{\mc{K}}(\mc{A}_{\chi})\right)^{{}^{(2)}\mathrm{A}},
\end{align}
where for the second equality we use the fact that the Hecke operators $\mathbb{T}_{\mc{K}}$ carries the sub-category $D^b(\tBun)_{\psi}$ to $D^b(\tBun \times X^0)_{\psi}$ for any $\psi \colon {}^{(2)}Z_G \to \Qlb^\times$ (recall from Lemma \ref{grouptheory} that $Z({}^{(2)}\mathrm{A})$ is equal to the double cover ${}^{(2)}Z_G \to Z_G$ of $Z_G= Z_G[2]$).

We have already observed that $(u_{\gamma} \times \id)^* \mathbb{T}_{\mc{K}}(\mc{A}_{\chi}) \cong (u_0 \times \id)^* \mathbb{T}_{\mc{K}}(\mc{A}_{\chi})$ is independent of $\gamma$; we conclude that
\[
\mathbb{T}_{\mc{K}}(\mc{A}_{\chi}) \cong \mc{A}_{\chi} \boxtimes \left( V_{\chi}^* \otimes (u_0 \times \id)^* \mathbb{T}_{\mc{K}}(\mc{A}_{\chi}) \right)^{{}^{(2)}\mathrm{A}},
\]
and we claim that $\mc{A}_{\chi}$ is a Hecke eigensheaf with `eigenvalue'
\begin{align}
\widetilde{\mc{E}}_{\chi} \colon \Sat_{\tG} &\to \Loc(X^0) \\
\mc{K}&\mapsto \left( V_{\chi}^* \otimes (u_0 \times \id)^* \mathbb{T}_{\mc{K}}(\mc{A}_{\chi}) \right)^{{}^{(2)}\mathrm{A}}.
\end{align}
That is, what remains to show is that $\tE_{\chi}(\mc{K})$ is in fact a local system, and that $\tE_{\chi}$ is a tensor functor satisfying the conditions of Definition \ref{eigensheafdefn}. This follows (by the monoidal property of the Hecke operators) by the same argument as \cite[\S 4.2]{heinloth-ngo-yun:kloosterman}, since we have seen that $\left( V_{\chi}^* \otimes (u_0 \times \id)^* \mathbb{T}_{\mc{K}}(\mc{A}_{\chi}) \right)^{{}^{(2)}\mathrm{A}}$ lies in perverse degree one.
\endproof
To summarize:
\begin{cor}\label{locsyslift}
Assume $k$ is as in Definition \ref{Vchi}. For every odd character $\chi \colon Z({}^{(2)}\mathrm{A})(\bar{k}) \to \Qlb^\times$, the object $\mc{A}_{\chi}$ of $D^b(\tBun)_{\mr{odd}}$ given by $\mc{A}_{\chi}|_{\tBun^{\gamma}}= j_{\gamma, !}(\mc{F}_{\chi})$ is a Hecke eigensheaf with eigen local system
\[
\widetilde{\mc{E}}_{\chi} \colon \Sat_{\tG} \to \Loc(X^0, \Ql'),
\]
giving rise by the Tannakian formalism to a monodromy representation (recall the notation from equation \ref{G1})
\[
\tilde{\rho}_{\chi} \colon \pi_1(X^0) \to \tG^\vee_1(\Ql').
\]
The restriction of $\widetilde{\mc{E}}_{\chi}$ to the full subcategory $\Sat_G \subset \Sat_{\tG}$ is naturally isomorphic to the eigen local system (there denoted $\mc{E}'_{\chi}$) of \cite[Theorem 4.2]{yun:exceptional}.

Moreover, if $\mc{K}= \IC_{\lambda}(m)$ is simple, then $\tE_{\chi}(\mc{K})$ is pure of weight $\langle 2 \rho, \lambda \rangle-2m$.
\end{cor}
\proof
We have established everything except the purity claim, which follows from the argument of Theorem \ref{eigensheaf}. Namely, equations (\ref{semiperverse}) and (\ref{dualsemiperverse}) imply that $\tE_{\chi}(\mc{K})$ is mixed of weights $\leq$ and $\geq$ $\langle 2\rho, \lambda, \rangle - 2m$ (by \cite{deligne:weil2}).
\endproof
Consequently, we have a commutative diagram
\[
\xymatrix{
& \tG^\vee_1(\Ql') \ar[d] \\
\pi_1(X^0) \ar[ur]^{\tilde{\rho}_{\chi}} \ar[r]^{\rho_{\chi}} & G^\vee_1(\Ql'),
}
\]
where $\rho_{\chi}$ is (of course these monodromy representations are only well-defined up to $\tG^\vee$, respectively $G^\vee$, conjugation) Yun's local system.
\section{The motives}
Having established the Hecke eigensheaf property, we can now describe the local systems $\tE_{\chi}(\mc{K})$ for all $\mc{K} \in \Sat_{\tG}$. We continue to assume $k$ is as in Definition \ref{Vchi}; in particular, the $k$-group scheme $Z({}^{(2)}\mathrm{A})$ is discrete. Let us fix a dominant weight $\lambda \in X^\bullet(\tT^\vee)= X_{\bullet}(\tT)$, and restrict to the case of $\mc{K}= \IC_{\lambda}$. In this case the sheaf $\mc{K}_{\tHk}$ is supported on a sub-stack $\tHk_{\leq \lambda}$, and the sheaf 
\[
\lefth^*(\mc{A}_{\chi}^0) \otimes \mc{K}_{\tHk}
\]
is supported on the locus of $(\mc{P}, \mc{P}', x, \iota)$ where $\mc{P} \in \tBun^0$ and $\mc{P}$ and $\mc{P}'$ are in relative position $\leq \lambda$, i.e. $\ev(\mc{P}, \mc{P}', \iota, x)$ lies in the $\leq \lambda$ strata of $\left[ (L^+ \tG \bs L \tG / L^+ \tG)/\Aut_{\mc{O}} \right]$. This forces $\mc{P}'$ to lie in the component $\tBun^{\nu \circ \lambda}$, where recall $\nu \colon \tG \to S$ is the multiplier character. It follows that to compute $\mathbb{T}_{\mc{K}}(\mc{A}_{\chi}^0)$ we can restrict $\lefth \colon \tHk_{\leq \lambda} \to \tBun$ to the preimage of $\tBun^0$, and thus consider instead the correspondence diagram
\[
\xymatrix{
& \tHk_{\leq \lambda}|_{\lefth^{-1}(\tBun^0)} \ar[ld]_{\lefth} \ar[dr]^{\righth} & \\
\tBun^0 & & \tBun^{\nu \circ \lambda} \times X^0.
}
\] 
In terms of this diagram, we find
\begin{equation}\label{eigencomp}
\mathbb{T}_{\mc{K}}(\mc{A}_{\chi}^0) \xrightarrow{\sim} \mc{A}_{\chi}^{\nu \circ \lambda} \boxtimes \tE_{\chi}(\mc{K}).
\end{equation}
Recall we are trying to describe $\tE_{\chi}(\mc{K})$. The argument is that of \cite[Lemma 4.3]{yun:exceptional}, except we have to keep track of the different connected components. Pulling back equation (\ref{eigencomp}) by $(u_{\nu \circ \lambda} \times \id)$, we obtain, just as in equations (\ref{maindiagram}) and (\ref{basechange}), a diagram with Cartesian squares
\begin{equation}\label{maindiagrambis}
\xymatrix{
& \widetilde{\mf{G}}^U_{\leq \lambda} \ar[ld] \ar[r]^{\upsilon_0} & \tGR^U_{u_{\nu \circ \lambda}, \leq \lambda} \ar[ld]_-{\omega^U_{u_{\nu \circ \lambda}}} \ar[r] & \tGR_{u_{\nu \circ \lambda}, \leq \lambda} \ar[ld]_-{\omega_{u_{\nu \circ \lambda}}} \ar[r] \ar[d] & \tHk_{\leq \lambda}|_{\lefth^{-1}(\tBun^0)} \ar[d]^{\righth} \\
\Spec k \ar[r]^{u_0} & [\twoK_0 \bs U] \ar[r] & \tBun^0 & X^0 \ar[r]^-{u_{\nu \circ \lambda} \times \id} & \tBun^{\nu \circ \lambda} \times X^0;
}
\end{equation}
and letting $\pi^U_{u_{\nu \circ \lambda}}$ denote the composite map $\tGR^U_{u_{\nu \circ \lambda}, \leq \lambda} \to X^0$, we obtain an identification
\begin{equation}\label{pulledback}
V_{\chi} \otimes \tE_{\chi}(\mc{K}) \cong (u_{\nu \circ \lambda} \times \id)^* \mathbb{T}_{\mc{K}}(\mc{A}_{\chi}^0) \cong \pi^U_{u_{\nu \circ \lambda}, !} \left( \omega^{U, *}_{u_{\nu \circ \lambda}}(\mc{F}_{\chi}) \otimes \mc{K}_{\tGR} \right).
\end{equation}
(We will write $\mc{K}_{\tGR}$ for the pull-back of $\mc{K}_{\tHk}$ to either of $\tGR_{u_{\nu \circ \lambda}, \leq \lambda}$ or $\tGR_{u_{\nu \circ \lambda}, \leq \lambda}^U$.) Also let 
\[
\pi_{\tmfG} \colon \tmfG \to X^0
\]
denote the corresponding projection. We now exploit the fact that $\widetilde{\mf{G}}^U_{\leq \lambda}$ carries a ${}^{(2)}\mathrm{A} \times {}^{(2)}\mathrm{A}$-action: for clarity the first copy, acting via the pull-back on the $\lefth$ (or as here, $\omega_{u_{\nu \circ \lambda}}$) projection, will be denoted ${}^{(2)}\mathrm{A}(1)$, and the second copy, acting via pull-back on the $\righth$ projection, will be denoted ${}^{(2)}\mathrm{A}(2)$. Decomposing the regular representation of ${}^{(2)}\mathrm{A}$, we obtain a ${}^{(2)} \mr{A}(1)$-equivariant isomorphism
\[
(\upsilon_{0, *} \Ql')_{\mr{odd}} \cong \bigoplus_{\substack{\chi \colon Z({}^{(2)}\mathrm{A}) \to \Qlb^\times \\ \text{$\chi$ is odd}}} V_{\chi}^* \otimes \omega_{u_{\nu \circ \lambda}}^{U, *} \mc{F}_{\chi}.
\]
Here ${}^{(2)}\mathrm{A}(1)$ acts on $V_{\chi}^*$. Since the isomorphism (\ref{pulledback}) is ${}^{(2)}\mathrm{A}(2)$-equivariant (acting on $V_{\chi}$ on the left-hand side, and on the right-hand side since $\mc{K}_{\tGR}$ is the pull-back of $\mc{K}_{\tHk}$), 
we obtain a ${}^{(2)}\mathrm{A} \times {}^{(2)}\mathrm{A}$-equivariant isomorphism
\begin{align}
\left( \pi_{\tmfG, !} \upsilon_0^* \mc{K}_{\tGR} \right)_{\mr{odd}} &\cong \left(\pi^U_{u_{\nu \circ \lambda}, !} \upsilon_{0, !}(\upsilon_0^* \mc{K}_{\tGR}) \right)_{\mr{odd}} \\
&\cong \bigoplus_{\substack{\chi \colon Z({}^{(2)}\mathrm{A}) \to \Qlb^\times \\ \text{$\chi$ is odd}}} V_{\chi}^* \otimes  \pi^U_{u_{\nu \circ \lambda, !}} (\omega^{U, *}_{u_{\nu \circ \lambda}} \mc{F}_{\chi} \otimes \mc{K}_{\tGR} ) \cong \bigoplus_{\substack{\chi \colon Z({}^{(2)}\mathrm{A}) \to \Qlb^\times \\ \text{$\chi$ is odd}}} V_{\chi}^* \otimes V_{\chi} \otimes \tE_{\chi}(\mc{K}).
\end{align}
Writing $\Ql'[{}^{(2)}\mathrm{A}]_{\chi}$ for the ${}^{(2)}\mathrm{A} \times {}^{(2)}\mathrm{A}$-equivariant local system on $\Spec k$ corresponding to the representation $V_{\chi}^* \otimes V_{\chi}$ of the group
\[
\left({}^{(2)}\mathrm{A}(\bar{k}) \times {}^{(2)}\mathrm{A}(\bar{k})\right) \rtimes \gal{k},\footnote{Equivalently, regarding $\Ql'[{}^{(2)}\mathrm{A}(\bar{k})]$ as ${}^{(2)}\mathrm{A}(\bar{k}) \times {}^{(2)}\mathrm{A}(\bar{k})$-module via $(a_1, a_2) \cdot a= a_1 a a_2^{-1}$, and extracting the constituent where $Z({}^{(2)}\mathrm{A}(2))$ acts by $\chi$.} 
\]
we summarize what we have shown (compare \cite[Lemma 4.3]{yun:exceptional}):
\begin{lemma}\label{calceigen}
There is a canonical isomorphism of ${}^{(2)}\mathrm{A} \times {}^{(2)}\mathrm{A}$-equivariant local systems on $X^0$
\[
\left( \pi_{\tmfG, !} \upsilon_0^* \mc{K}_{\tGR} \right)_{\mr{odd}} \cong \bigoplus_{\substack{\chi \colon Z({}^{(2)}\mathrm{A}) \to \Qlb^\times \\ \text{$\chi$ is odd}}} \Ql'[{}^{(2)}\mathrm{A}]_{\chi} \otimes \tE_{\chi}(\mc{K}).
\]
In particular, the left-hand side is a local system.
\end{lemma}
It is explained in \cite[\S 3.3.4]{yun:exceptional} how to take the invariants of an equivariant perverse sheaf under a (not necessarily discrete) finite group scheme. Applying this:
\begin{cor}
For all $\mc{K} \in \Sat_{\tG}$ and all odd $\chi \colon Z({}^{(2)}\mathrm{A}) \to \Qlb^{\times}$, there is an isomorphism of local systems on $X^0$:
\[
\tE_{\chi}(\mc{K}) \cong \left(\Ql'[{}^{(2)}\mathrm{A}]_{\chi}^* \otimes (\pi_{\widetilde{\mf{G}}^U_{\leq \lambda}, !} \upsilon_0^* \mc{K}_{\tGR})_{\mr{odd}} \right)^{{}^{(2)}\mathrm{A} \times {}^{(2)}\mathrm{A}}.
\]
\end{cor}
\section{The case of minuscule weights}\label{minuscule}
We now want to make this description of $\tE_{\chi}(\mc{K})$ explicit. Our ultimate goal is the following:
\begin{thm}\label{motivated}
Let $k$ be $\Q$ or $\Q(\sqrt{-1})$ according as $G$ is of type $D_{4n}, G_2, E_8$ or $A_1, D_{4n+2}, E_7$. Consider any odd $\chi \colon Z({}^{(2)}\mathrm{A}) \to \Qlb^\times$ and any $\mc{K} \in \Sat_{\tG}$. Let $F$ be any number field containing $k$. Then for any point $t \in X^0(F)$, the specialization
\[
\tilde{\rho}_{\chi, \mc{K}, t} \colon \pi_1(\Spec F) \xrightarrow{t} \pi_1(X^0) \xrightarrow{\tilde{\rho}_{\chi}} \tG^\vee_1(\Ql') \to \mr{GL}\left(H^*_w (\mc{K})\right)
\]
(the representation of $\tG^\vee_1$ is that induced by $\mc{K}$ under the Satake isomorphism, as in \S \ref{satake}) is, as $\gal{F}$-representation, isomorphic to the $\Ql'$-realization of an object of $\mc{M}_{F, \Q'}$.
\end{thm} 
The case of $\mc{K}$ corresponding to a quasi-minuscule weight is considered in \cite[\S 4.3]{yun:exceptional}. Although our discussion is valid for any $\tG$ as in \S \ref{setting}, there are certain cases in which it is uninteresting: for instance, if $G= \mr{SL}_2$, we gain nothing by taking $\tG= \mr{SL}_2 \times \mathbb{G}_m$; however, by taking $\tG= \mr{GL}_2$, we gain the representations of $\mr{SL}_2= G^\vee_{\mr{sc}}$ (the simply-connected cover of $G^\vee$), and it is these new representations that will be of interest, just as in the classical setting the Kuga-Satake abelian variety is found via the spin representation of $\mr{Spin}_{21}$, while the motive of the $K3$ arises from the standard 21-dimensional representation. 

To show that $\tilde{\rho}_{\chi, \mc{K}, t}$ arises from an object of $\mc{M}_{F, \Q'}$ demands a significant digression into understanding intersection cohomology of varieties with arbitrarily bad singularities. A good first approximation to understanding the motivic nature of $\tilde{\rho}_{\chi, t} \colon \gal{F} \to \tG^\vee_1(\Ql')$ is to verify this after composition with a single \textit{faithful} finite-dimensional representation of $\tG^\vee_1$ (i.e., to show $\tilde{\rho}_{\chi, t}$ is weakly motivic in the sense of Definition \ref{genKSdefn}). That is what we will do in this section.

First, we make a robust choice of $\tG$, such that $\tG^\vee$ has representations restricting to each of the minuscule representations of $G^\vee_{\mr{sc}}$. For instance, we can take:
\begin{itemize}
\item ($A_1$) $\tG= \mr{GL}_2$;
\item ($E_7$) Let $c$ denote the non-trivial element of $Z_G= \mu_2$. Then take $\tG= (G \times \mathbb{G}_m)/\langle (c, -1) \rangle$;
\item ($D_n$, $n$ even) Let $c$ and $z$ be generators of $Z_G \cong \mu_2 \times \mu_2$. Then take 
\[
\tG= (G \times \mathbb{G}_m \times \mathbb{G}_m)/\langle (c, -1, 1), (z, 1, -1) \rangle.
\]
\end{itemize} 
Each minuscule representation of $G^\vee_{\mr{sc}}$\footnote{These are, in the three cases: the standard representation of $\mr{SL}_2$, the 56-dimensional representation of $E_7$, and the standard ($2n$-dimensional) and two half-spin representations of $\mr{Spin}_{2n}$.} extends to an irreducible representation of $\tG^\vee$, and then to an irreducible representation of $\tG^\vee_1$. Taking a direct sum, we obtain a faithful family of representations
\begin{equation}\label{faithfulfamily}
r_{\mr{min}} \colon \tG^\vee_1 \to \mr{GL}(V_{\mr{min}}),
\end{equation}
and set ourselves the goal of showing that each $r_{\mr{min}} \circ \tilde{\rho}_{\chi, t}$ is motivated. The full force of Theorem \ref{motivated} is considerably deeper (it is new even in Yun's original setting), so in the present section we will only treat the case of these minuscule weights, which has the added advantage that the relevant geometry--of the corresponding affine Schubert varieties--is especially simple. 

We begin, however, with some generalities: continue to let $\mc{K} \in \Sat_{\tG}$ be any irreducible object of the form $\mc{K}= \IC_{\lambda}$, $\lambda \in X^\bullet(\tT^\vee)$ (the discussion will apply equally well to $\mc{K}$ of the form $\IC_{\lambda}(m)$, but we take $m=0$ to simplify the notation). What we denoted above by $\mc{K}_{\tGR}[1]$ is the intersection complex of $\tGR_{u_{\nu \circ \lambda}, \leq \lambda}^U$ (or the same before restricting to $U$). Since the map 
\[
\upsilon_0 \colon \widetilde{\mf{G}}^U_{\leq \lambda} \to \tGR_{u_{\nu \circ \lambda}, \leq \lambda}
\] 
is \'{e}tale, $\upsilon_0^* \mc{K}_{\tGR}[1]$ is again the intersection complex of $\widetilde{\mf{G}}^U_{\leq \lambda}$. Recall that the stratification of the affine Grassmannian induces one for the Beilinson-Drinfeld Grassmannian:
\[
\tGR_{u_{\nu \circ \lambda}, \leq \lambda}= \coprod_{\substack{\mu \leq \lambda \\ \text{$\mu$ dominant}}} \tGR_{u_{\nu \circ \lambda}, \mu}.
\]
The terms on the right-hand side are defined by replacing $\tHk_{\leq \lambda}$ by $\tHk_{\mu}$ in diagram (\ref{maindiagrambis}). Note that $\nu \circ \mu= \nu \circ \lambda$ since $\lambda- \mu \in X_\bullet(T)$ lies in the co-root lattice of $G$. The dense open locus $\tGR_{u_{\nu \circ \lambda}, \lambda}$ is smooth over $X^0$: fiber-wise it is the smooth stratum $\Gr_{\tG, \lambda}$ of $\Gr_{\tG, \leq \lambda}$. We write $\widetilde{\mf{G}}^U_{\lambda}$ and $\widetilde{\mf{G}}^U_{< \lambda}$ for the preimages in $\widetilde{\mf{G}}^U_{\leq \lambda}$ of $\tGR_{u_{\nu \circ \lambda}, \lambda}$ and $\coprod_{\mu < \lambda} \tGR_{u_{\nu \circ \lambda}, \mu}$.

Taking the $t$-fiber ($t \in X^0(F)$) of the isomorphism in Lemma \ref{calceigen}, we obtain a (quasi-)isomorphism
\begin{equation}\label{IHc}
\IH_c(\widetilde{\mf{G}}^U_{\leq \lambda, t})_{\mr{odd}} \cong \bigoplus_{\text{$\chi$ odd}} \Ql'[{}^{(2)}\mathrm{A}]_{\chi} \otimes \tilde{\rho}_{\chi, \mc{K}, t}.
\end{equation}
Let us explain the notation. For any irreducible variety $Y$ over a field $F$, the intersection complex $\IC_Y$ is a perverse sheaf in cohomological degrees $[-\dim Y, 0]$. It is pure of weight $\dim Y$. We denote by $\IH_c(Y)$ the complex $R\Gamma_c(\IC_Y)$ on $\Spec F$: it lies in cohomological degrees $[-\dim Y, \dim Y]$, and is pure of weights $\leq \dim Y$. As usual, we then define the compactly-supported intersection cohomology $\IH_c^i(Y_{\overline{F}})$ (a $\gal{F}$-representation) as $H^{i-\dim Y}(\IH_c(Y))$ (note the degree shift). We also observe that while (compactly-supported) intersection cohomology is not in general functorial for (proper) morphisms of varieties, it is for (proper) \'{e}tale morphisms: $\IC_{\widetilde{\mf{G}}^U_{\leq \lambda, t}}$ is still ${}^{(2)}\mathrm{A} \times {}^{(2)} A$-equivariant, as is the isomorphism (\ref{IHc}). Now as in \cite[\S 4.3.2]{yun:exceptional}, we let $e_{\chi} \in \Q'[{}^{(2)}\mathrm{A}(\bar{k}) \times {}^{(2)}\mathrm{A}(\bar{k})]^{\gal{k}}$ be the idempotent whose action on the ${}^{(2)}\mathrm{A}(\bar{k})\times {}^{(2)}\mathrm{A}(\bar{k})$-module $\Q'[{}^{(2)}\mathrm{A}(\bar{k})]$ projects to the component $\Q'[{}^{(2)}\mathrm{A}]_{\chi}$ and then onto the line spanned by $\id \in \End(V_{\chi})$ (a direct factor of the representation $\Q'[{}^{(2)}\mathrm{A}]_{\chi}$ after restricting to the diagonal copy ${}^{(2)}\mathrm{A}(\bar{k}) \into {}^{(2)}\mathrm{A}(\bar{k}) \times {}^{(2)}\mathrm{A}(\bar{k})$). Explicitly,
\[
e_{\chi}= \frac{1}{|{}^{(2)} A(\bar{k}) \times {}^{(2)} A(\bar{k})|} \sum_{(a_1, a_2)} \theta_{\chi}(a_1 a_2^{-1}) (a_1, a_2),
\]
where $\theta_{\chi}$ denotes the character of the ${}^{(2)} A(\bar{k})$-representation $V_{\chi}$.
\begin{prop}\label{explicit}
Let $\mc{K}= \IC_{\lambda} \in \Sat_{\tG}$, let $\chi \colon Z({}^{(2)}\mathrm{A}) \to \Qlb^\times$ be odd, and let $t \in X^0(F)$ for any number field $F$ containing $k$. Then
\[
\tilde{\rho}_{\chi, \mc{K}, t} \cong \Gr^W_{\langle 2\rho, \lambda \rangle} \left( e_{\chi} \IH^{\langle 2 \rho, \lambda \rangle}_c(\widetilde{\mf{G}}^U_{\leq \lambda, t})\right).
\]
\end{prop}
\proof
Apply $e_{\chi}$ to equation (\ref{IHc}), noting that the right-hand side is concentrated in degree zero. Since we have seen that $\tE_{\chi}(\IC_{\lambda})$ is pure of weight $\langle 2 \rho, \lambda \rangle$, the claim is immediate.
\endproof
Proposition \ref{explicit} reduces Theorem \ref{motivated} to a special case of the following general theorem:
\begin{thm}\label{ICmotivated}
Let $k$ be a finitely-generated field of characteristic zero, and let $Y/k$ be a quasi-projective variety. Then for all $i, r \in \Z$, $\Gr^W_i(\IH^r_c(Y))$ is as $\gal{k}$-representation isomorphic to the $\ell$-adic realization of an object of $\mc{M}_k$. 

Next consider the case in which $Y$ is acted on by a finite $k$-group scheme $\Gamma$. Let $e \in \Qb[\Gamma(\bar{k})]^{\gal{k}}$ be an idempotent. Fix an embedding $\Qb \into \Qlb$. Then for all $i, r \in \Z$, $\Gr^W_i(e \IH^r_c(Y, \Qlb))$ is as $\gal{k}$-representation isomorphic to the $\Qb \into \Qlb$-realization of an object of $\mc{M}_{k, \Qb}$.
\end{thm}
\begin{rmk}
See \S \ref{overview} for what is meant by the weight-gradings $\Gr^W_{\bullet}$. Note that in the application we only need the case $i=r$.
\end{rmk}
Theorem \ref{ICmotivated} will be proven in Corollary \ref{ICmotpf}. For the remainder of this section, we content ourselves with showing that $\tilde{\rho}_{\chi, t}$ is weakly motivic. Thus, it suffices to assume that $\lambda$ restricts to a minuscule weight of $G^\vee_{\mr{sc}}$. In this case, $\Gr_{\tG, \leq \lambda}= \Gr_{\tG, \lambda}$ has non-singular reduced part, so that
\[
\tilde{\rho}_{\chi, \IC_{\lambda}, t} \cong \Gr^W_{\langle 2 \rho, \lambda \rangle} \left(e_{\chi} H^{\langle 2 \rho, \lambda \rangle }_c(\widetilde{\mf{G}}^U_{\lambda, t}) \right).
\]
That the right-hand side is isomorphic to the $\ell$-adic realization of an object of $\mc{M}_{F, \Q'}$ follows from the standard description (originating in \cite{deligne:ht1}) of the weight filtration on the cohomology of a smooth variety, via the Leray spectral sequence for its inclusion into a smooth compactification with boundary given by a smooth normal crossings divisor. See \cite[\S 4.3.1]{yun:exceptional} or \cite{patrikis-taylor:irr} (the discussion between Remark 2.6 and Lemma 2.7) for this equivariant version. We conclude:
\begin{cor}
For all choices of $\tG$ as in \S \ref{setting}, there exists a faithful finite-dimensional representation $r \colon \tG^\vee_1 \into \mr{GL}(V_r)$ such that for all number field specializations $\Spec F \xrightarrow{t} X^0$ with $F$ satisfying condition (\ref{khassqrt}), 
\[
r \circ \tilde{\rho}_{\chi, t} \colon \gal{F} \to \mr{GL}(V_r \otimes \Ql')
\]
is isomorphic to the $\Ql'$-realization of an object of $\mc{M}_{F, \Q'}$. For all $G$, we may choose $\tG$ and $r$ such that $r|_{G^\vee_{\mr{sc}}}$ is isomorphic to the direct sum of all the minuscule representations of $G^\vee_{\mr{sc}}$.

In particular, the lifts $\tilde{\rho}_{\chi, t}$ of Yun's $\rho_{\chi, t}$ satisfy the generalized Kuga-Satake property of Definition \ref{genKSdefn}.
\end{cor}
\section{Intersection cohomology is motivated}
\subsection{Overview}\label{overview}
In the remaining sections, which are logically independent of the rest of the paper, we prove Theorem \ref{ICmotivated}. Let $k$ be a field of characteristic zero, and fix an algebraic closure $\bar{k}$ of $k$. As usual let $\gal{k}= \Gal(\bar{k}/k)$. Let $Y/k$ be any quasi-projective variety. If $Y$ is irreducible of dimension $d_Y$, we can form the $\ell$-adic intersection cohomology groups 
\[
\IH^{r+d_Y}(Y)= H^r(Y_{\bar{k}}, \IC_Y|_{Y_{\bar{k}}}), 
\]
as well as their analogues with compact supports, $\IH^{r+d_Y}_c(Y)$. If $Y$ is reducible, the definitions need a little more care, working component by component; see \cite[\S 4.6]{decataldo:pervfil} for an explanation. The intersection complex $\IC_{Y_{\bar{k}}}$ is $\gal{k}$-equivariant, so $\gal{k}$ acts on $\IH^*(Y)$ and $\IH^*_c(Y)$. Since we do not assume $Y$ is projective, these $\gal{k}$-representations are not pure; in particular, Theorem \ref{ICmotivated} cannot hold for the groups $\IH^*_c(Y)$ themselves. Thus we first need to make sense of the weight filtration on $\IH^*_c(Y)$, in order even to speak of the $\gal{k}$-representations $\Gr^W_{\bullet} \IH^*_c(Y)$.

There are two basic templates, one `sheaf-theoretic' and one `geometric,' for endowing the cohomology of a variety with a weight filtration. The models for the former approach are \cite{deligne:weil2} and \cite{bbd}; the models for the latter are \cite{deligne:ht2} and \cite{deligne:ht3}. The latter approach typically depends on having resolution of singularities over the field $k$, and is consequently restricted to characteristic zero; but when available, it yields more robust, because more `motivic,' results. Thus we will explain, at least for $k$ finitely-generated over $\Q$, how to give an \textit{a priori} `sheaf-theoretic' sense to $\Gr^W_{\bullet} \IH^*(Y)$, but then our main aim will be to give a `geometric' construction, as part of the proof of Theorem \ref{ICmotivated}, that recovers the sheaf-theoretic definition of the $\gal{k}$-representations $\Gr^W_{\bullet} \IH^*_c(Y)$. Let us begin then by recalling the sheaf-theoretic construction of a weight filtration on $\IH^*_c(Y)$.

Since we work with $k$ of characteristic zero, the basic case of positive characteristic addressed in \cite{deligne:weil2} and \cite{bbd} is not sufficient. But the results of those papers have been extended in a form suitable for our purposes, and indeed much more generally than we require, in \cite{huber:mixedperverse} and \cite{morel:weights}.\footnote{The basic notions of horizontal sheaf, perverse t-structure on the `derived' category of horizontal sheaves, and weights for horizontal sheaves are developed in, respectively, \S 1, 2, and 3 of \cite{huber:mixedperverse}. Morel's paper builds on these foundations, generalizing the results of \cite{huber:mixedperverse} to any finitely-generated $k$, and establishing a sort of six operations functoriality for complexes having weight filtrations.} Namely, the intersection complex $\IC_Y$ is a horizontal, pure perverse sheaf in the sense of \cite[\S 2]{morel:weights}, and \cite[Th\'{e}or\`{e}me 3.2, Proposition 6.1]{morel:weights} implies that $\IH_c^*(Y)$ (likewise $\IH^*(Y)$) carries a unique weight filtration $W_{\bullet}$. In particular, this means that each $\Gr^W_r \IH^*_c(Y)$ is pure of weight $r$ in the following sense: the underlying lisse sheaf on $\Spec k$ arises by base-change from a lisse sheaf $\mc{G}$ on some smooth sub-algebra $A \subset k$, of finite-type over $\Z$, and with $\Frac(A)=k$; and for all specializations at closed points $x$ of $\Spec A$, $x^* \mc{G}$ is pure of weight $r$ in the usual finite field sense. This characterizing property will hold for the output of our geometric construction; this is verified step-by-step as the construction proceeds.   

We now outline the approach to Theorem \ref{ICmotivated}. By Poincar\'{e} duality for intersection cohomology (which is $\gal{k}$-equivariant), we may restrict to the case of $\IH^*(Y)$. First, we remark that the basic difficulty, and interest, of this problem is that both intersection cohomology and weight filtrations are \textit{a priori} `sheaf-theoretically' defined. The theorem shows that these sheafy constructions can in fact be realized just by playing with the cohomology of smooth projective varieties. There are two, essentially orthogonal, special cases of this problem:
\begin{itemize}
\item $Y$ may be smooth but non-projective. In this case, $\IH^r(Y)= H^r(Y_{\bar{k}}, \Ql)$, and the result follows from the geometric approach of \cite{deligne:ht2}; namely, if $\overline{Y}$ is a smooth compactification of $Y$ with $\overline{Y} \setminus Y$ equal to a union of smooth divisors $D_{\alpha}$ with normal crossings, then the ($E_3$-degenerate) Leray spectral sequence for the inclusion $Y \subset \overline{Y}$ yields a description of $\Gr^W_{\bullet} H^r(Y_{\bar{k}}, \Ql)$ in terms of the divisors $D_{\alpha}$ and their various (smooth, projective) intersections; see Theorem \ref{hfunctor}, part 3, below, for a slight rephrasing.
\item $Y$ may be projective but singular. In this case, the result, when $k$ is algebraically closed,\footnote{In this case one should work not just with $\ell$-adic cohomology but also with (compatible) Betti and de Rham realizations, in order for the assertion to have any content.} has been proven by de Cataldo and Migliorini. We briefly describe the two crucial geometric inputs (assume for this informal description that $k$ is algebraically closed). Let $f \colon X \to Y$ be a resolution of singularities. Roughly speaking, $\IH^*(Y)$ occurs as a `main term' in $H^*(X, \Ql)= H^*(Y, f_* \Ql)$ corresponding (via the decomposition theorem) to the summand of $f_* \Ql$ in perverse degree $\dim X$ and supported along the open dense stratum (the non-singular locus) $Y^0$ of $Y$. The first key result is that the perverse (Leray) filtration on $H^*(Y, f_* \Ql)$ admits (\cite[Theorem 3.3.5]{decataldo:pervfil}) a remarkable geometric description in terms of a suitably generic `flag filtration.' The second is that the factor of $f_* \Ql$ supported along $Y^0$ can, at least in cohomology, also be extracted `geometrically': this follows from the novel approach to the decomposition theorem pioneered by de Cataldo-Migliorini in a series of papers (see \cite[\S 1.3.3]{decataldo-migliorini:decomproj} for a precise statement).
\end{itemize}
Our task is to fuse these two approaches, and to get everything to work over an arbitrary (not algebraically closed) field $k$ of characteristic zero. The chief obstruction to getting the relevant arguments of \cite{decataldo-migliorini:decomproj} to work over any $k$ is that the `generic flags' mentioned above would need to be defined $k$-rationally. This it turns out is not so hard to achieve, using Bertini's theorem over $k$ and, crucially, the fact that flag varieties are rational, so that any Zariski open over $k$ necessarily has $k$-points.

Rather more complicated is integrating the approaches of \cite{deligne:ht2} and \cite{decataldo-migliorini:decomproj} in order to prove Theorem \ref{ICmotivated} for any quasi-projective $Y$. The basic difficulty is that, since motivated motives are only defined in the pure case, the argument (resting on \cite{deligne:ht2}) in the smooth case is not obviously `functorial in $Y$.' Fortunately, it can be upgraded to one that is, using the results of \cite{guillen-navarro:criterion} on the existence of `weight complexes' of motivated motives whose cohomology computes $\Gr^W_{\bullet} H^*(Y)$ for any $k$-variety $Y$. We will also use a version for cohomology with compact supports: the latter, due independently to Gillet-Soul\'{e} (\cite{gillet-soule:weights}) and Guillen-Navarro, is somewhat simpler, but not suited for describing the perverse Leray filtration as in \cite{decataldo:pervfil}, even for cohomology with compact supports. It is crucial, however, that we exploit both theories: the inductive construction of the support decomposition as in \cite[Proposition 2.2.1]{decataldo-migliorini:decomproj} requires having motivated versions both of pull-back in $H^*$ and pull-back for proper morphisms in $H^*_c$ (note that these two kinds of pull-backs are not related by Poincar\'{e} duality; one cannot be formally reduced to the other). Once this setup is in place, however, the arguments of \cite{decataldo-migliorini:decomproj} go through \textit{mutatis mutandis}. We consequently establish stronger results on finding `motivated' splittings of the perverse Leray filtration, and a motivated support decomposition, closely in parallel to the main results of \cite{decataldo-migliorini:decomproj}: see Theorem \ref{motivateddecomposition} and Corollary \ref{splitting}, which should be regarded as the main results of this half of the paper. 
\begin{notation}
Except where we explicitly allow more general fields, from now on $k$ will be a finitely-generated field extension of $\Q$. Whenever we speak of the weight-grading $\Gr^W_{\bullet}$ on various cohomology groups of a variety over $k$, the grading is unique, and can be shown to exist by \cite[Th\'{e}or\`{e}me 3.2, Proposition 6.1]{morel:weights}. As before, $\mc{M}_k$ denotes Andr\'{e}'s category of motives for motivated cycles over $k$ (with $\Q$-coefficients). For a smooth projective variety $X$ over $k$, we write $H(X)$ for the canonical object of $\mc{M}_k$ associated to $X$. Finally, given a map of varieties $f \colon X \to Y$, we always mean the derived functors when we write $f_*$, $f_!$, etc.
\end{notation}
\subsection{Weight-graded motivated motives associated to smooth varieties}
Here is the theorem of Guill\'{e}n and Navarro Aznar, specialized to the precise statement we require:\footnote{They prove something stronger, with Chow motives in place of motivated motives.}
\begin{thm}[see Th\'{e}or\`{e}me 5.10, \cite{guillen-navarro:criterion}]\label{hfunctor}
Let $k$ be a field of characteristic zero, and let $\mr{Sch}/k$ denote the category of finite-type separated $k$-schemes. Then there exists a contravariant functor
\begin{equation}
h \colon \mr{Sch}/k \to K^b(\mc{M}_k),
\end{equation}
valued in the homotopy category of bounded complexes in $\mc{M}_k$, such that
\begin{enumerate}
\item If $X$ is a smooth projective $k$-scheme, then $h(X)$ is naturally isomorphic to the canonical motivated motive $H(X)$ associated to $X$.
\item If $X$ is a smooth projective $k$-scheme, and $D= \cup_{\alpha=1}^t D_{\alpha}$ is a normal crossings divisor equal to the union of smooth divisors $D_{\alpha}$, we can form a cubical diagram of smooth projective varieties
\[
S_{\bullet}(D) \to X,
\]
where for every non-empty subset $\Sigma \subset \{1, \ldots, t\}$, $S_{\Sigma}(D)$ is the (smooth) intersection $D_{\Sigma}= \cap_{\alpha \in \Sigma} D_{\alpha}$, with the obvious inclusion maps $S_{\Sigma}(D) \to S_{\Sigma'}(D)$ whenever $\Sigma' \subset \Sigma$. Using the \textbf{covariant} functoriality arising from Gysin maps, we can then associate a cubical diagram $h_*(S_{\bullet}(D) \to X)$ in $\mc{M}_k$; to be precise, $h_*(S_{\Sigma}(D))$ is the object of $\mc{M}_k$
\[
h(D_{\Sigma})(\dim D_{\Sigma}),
\]
with the Gysin maps $h_*(S_{\Sigma}(D)) \to h_*(S_{\Sigma'}(D))$ whenever $\Sigma' \subset \Sigma$. Then $h(X \setminus D)$ is isomorphic to the simple complex associated to this cubical diagram (see the proof for what this means):
\begin{equation}\label{opensmoothcalc}
h(X \setminus D) \cong \mbf{s}\left(h_*(S_{\bullet}(D) \to X)\right)(-\dim X).
\end{equation} 
\item In particular, $h(X \setminus D)$ is a complex whose degree $r$ homology\footnote{We use homological conventions here.} $H_r(h(X \setminus D))$ is an object of $\mc{M}_k$ whose $\ell$-adic realization is given by (for $k$ finitely-generated over $\Q$)
\[
H_{\ell}\left( H_r(h(X \setminus D)) \right) \cong \bigoplus_q \Gr^W_{q+r} H^q((X\setminus D)_{\bar{k}}, \Ql).
\]
\end{enumerate}
\end{thm}
\proof
Except for the third assertion, this is all explicitly in \cite[Th\'{e}or\`{e}me 5.10]{guillen-navarro:criterion}. The remaining claim follows from the usual (\cite{deligne:ht2}) description of the weight gradeds for $H^*((X \setminus D)_{\bar{k}}, \Ql)$: ignoring for notational convenience the Tate twists, the degree $r$ term $h(X \setminus D)_r$ (to be precise, after the identification of equation (\ref{opensmoothcalc})) is $\oplus_{|\Sigma|=r|} h(D_{\Sigma})$, with the boundary map $h(X \setminus D)_r \to h(X \setminus D)_{r-1}$ given by an alternating sum of Gysin maps. The $\ell$-adic realization of this complex can be identified (up to a sign in the boundary maps, at least: see \cite[1.8 Proposition]{guillen-navarro:locinvcyc}) with the complex
\[
\cdots \to K_r= \bigoplus_q E_1^{-r, q+r} \xrightarrow{\oplus_q d_1^{-r, q+r}} K_{r-1}= \bigoplus_q E_1^{-r+1, q+r} \to \cdots
\]
built out of the $E_1$ terms of the (weight) spectral sequence of the filtered complex (b\^{e}te filtration) 
\[
E_1^{-r, q+r}= H^q(X_{\bar{k}}, \Gr^W_r j_* \Ql) \implies H^q(X_{\bar{k}}, j_* \Ql)= H^q( (X \setminus D)_{\bar{k}}, \Ql).
\]
This spectral sequence degenerates at the $E_2$ page (by the yoga of weights), and its $E_2$ terms then give the weight gradeds of $H^q((X \setminus D)_{\bar{k}}, \Ql)$; part (3) of the Theorem follows.
\endproof
This is not a full description of the result of Guill\'{e}n and Navarro Aznar, but it contains the two points of interest for us: the explicit description of the objects $H_r(h(X \setminus D))$, and in particular their connection with the weight filtration on $H^*((X \setminus D)_{\bar{k}}, \Ql)$; and, crucially, the fact that $h$ is functorial. In particular, for any morphism $\phi \colon U \to V$ in $\mr{Sch}/k$, we get, for all $r$, morphisms $H_r(h(V)) \to H_r(h(U))$ in $\mc{M}_k$.

Here is the compact supports version:
\begin{thm}[Theorem 2 of \cite{gillet-soule:weights}, or Th\'{e}or\`{e}me 5.2 \cite{guillen-navarro:criterion}]\label{Wfunctor}
Let $k$ be a field of characteristic zero, and let $\mr{Sch}_c/k$ denote the category of separated finite-type $k$-schemes with morphisms given by proper maps. Then there exists a contravariant functor
\begin{equation}
W \colon \mr{Sch}_c/k \to K^b(\mc{M}_k)
\end{equation}
such that
\begin{enumerate}
\item If $X$ is a smooth projective $k$-scheme, then $W(X)$ is naturally isomorphic to the canonical motivated motive $H(X)$ associated to $X$.
\item If $X$ is a smooth projective $k$-scheme, and $D= \cup_{\alpha=1}^t D_{\alpha}$ is a normal crossings divisor equal to the union of smooth divisors $D_{\alpha}$, then $W(X \setminus D)$ is isomorphic to the simple complex (we now use cohomological conventions and normalize $W(X \setminus D)$ to live in cohomological degrees $[0, t]$)
\begin{equation}
H(X) \to \oplus_{\alpha} H(D_{\alpha}) \to \cdots \to \oplus_{|\Sigma|=s} H(D_{\Sigma}) \to \cdots,
\end{equation}
with coboundaries given by an alternating sum of restriction maps $H(D_{\Sigma'}) \to H(D_{\Sigma})$ whenever $\Sigma' \subset \Sigma$. (See \cite[Proposition 3]{gillet-soule:weights}.)
\item In particular, $W(X \setminus D)$ is a complex whose degree $s$ cohomology $H^s(W(X \setminus D))$ is an object of $\mc{M}_k$ whose $\ell$-adic realization is given by (for $k$ finitely-generated over $\Q$)
\[
H_{\ell} \left( H^s(W(X \setminus D))\right) \cong \bigoplus_p \Gr^W_p H^{p+s}_c((X \setminus D)_{\bar{k}}, \Ql).
\]
\end{enumerate}
\end{thm}
In the setting of parts 2 and 3 of Theorems \ref{hfunctor} and \ref{Wfunctor}, let $U= X \setminus D$. Poincar\'{e} duality for $U$ descends to a duality relation in $\mc{M}_k$ between the cohomologies of the complexes $h(U)$ and $W(U)$. Before stating it, we introduce a little more notation:
\begin{defn}
Let $H^q_r(h(U))$ be the canonical summand of $H_r(h(U))$ in $\mc{M}_k$ of weight $q+r$. Let $W^p(U)$ be the canonical complex of weight $p$ summands of the terms of $W(U)$, and let $H^s(W^p(U))$ be the degree $s$ cohomology. 
\end{defn}
\begin{rmk}
The object $H^q_r(h(U))$ of $\mc{M}_k$ has $\ell$-adic realization $\Gr^W_{q+r} H^q(U_{\bar{k}}, \Ql)$. The object $H^s(W^p(U))$ of $\mc{M}_k$ has $\ell$-adic realization $\Gr^W_p H^{p+s}_c(U_{\bar{k}}, \Ql)$.
\end{rmk}
\begin{lemma}\label{motivatedopenPD}
Let $U= X \setminus D$ as above, and assume $U$ is equidimensional of dimension $d$. Then there is a canonical isomorphism in $\mc{M}_k$
\begin{equation}
H^q_r(h(U))^\vee \cong H^r(W^{2d-q-r}(U))(d)
\end{equation}
\end{lemma}
\proof
Poincar\'{e} duality for each $D_{\Sigma}$ induces a perfect duality between $h(U)_s$ and $W(U)^s$ (the degree $s$ terms of each complex) for all $s$. The Gysin maps $H(D_{\Sigma'})(d-|\Sigma'|) \to H(D_{\Sigma})(d-|\Sigma|)$ are Poincar\'{e} dual to the pull-back maps $H(D_{\Sigma}) \to H(D_{\Sigma'})$ for all $\Sigma \subset \Sigma'$, and we can deduce perfect dualities (in $\mc{M}_k$ )
\[
H_s(h(U))^\vee \cong H^s(W(U))(d).
\]
The result follows from decomposing these dualities into each of their graded components.
\endproof
\section{The perverse Leray filtration}
\subsection{Relation to flag filtrations}
In this section we recall the beautiful and fundamental result of de Cataldo-Migliorini and de Cataldo that describes the perverse Leray filtration for a map of varieties in terms of a certain flag filtration: see \cite{decataldo-migliorini:pervfil} and, for the specific result we use, \cite[Theorem 3.3.5]{decataldo:pervfil}. These results are worked out over an algebraically closed field of characteristic zero, and so our first aim in this section is to check the analogue for any $k \supset \Q$. 

We first recall the Jouanolou trick.
\begin{defn}
Let $Y$ be a variety over $k$. An \textit{affinement} of $Y$ is a map $\mc{Y} \xrightarrow{p} Y$ in $\mr{Sch}/k$ with $\mc{Y}$ an affine $k$-scheme, such that $p$ is a torsor for some vector bundle on $Y$. 
\end{defn}
\begin{prop}[Lemme 1.5 of \cite{jouanolou:trick}]
Suppose $Y \in \mr{Sch}/k$ is quasi-projective. Then an affinement of $Y$ exists.
\end{prop}
Jouanolou's result in fact holds for arbitrary quasi-projective schemes, but we are only interested in the case of varieties over $k$. 

Now let $f \colon X \to Y$ be a morphism of $k$-varieties. It induces the (increasing) perverse Leray filtration on $H^*(X_{\bar{k}}, \Ql)$ via
\begin{equation}\label{defpervfil}
\mc{P}^f_j \left(H^*(X_{\bar{k}}, \Ql) \right)= \im \left( H^*(Y_{\bar{k}}, {}^p \tau_{\leq j} f_* \Ql) \to H^*(Y_{\bar{k}}, f_* \Ql) \right) \subseteq H^*(Y_{\bar{k}}, f_* \Ql)=H^*(X_{\bar{k}}, \Ql).
\end{equation}
Here ${}^p \tau_{\leq j}$ denotes perverse truncation.\footnote{It is defined for the complex $f_* \Ql$ on $Y$ itself, and we omit the base-change to $\bar{k}$ in the notation of these cohomology groups.} We make the analogous definition of the perverse Leray filtration on $H^*_c(X_{\bar{k}}, \Ql)$, replacing $f_* \Ql$ by $f_! \Ql$ (the only case of interest to us will be when $f$ is proper, so $f_*=f_!$).
\begin{thm}[Theorem 3.3.5 of \cite{decataldo:pervfil}]\label{pervleraythm}
Assume $k= \bar{k}$ is an algebraically closed field of characteristic zero. Let $f \colon X \to Y$ be a morphism in $\mr{Sch}/k$ with $Y$ quasi-projective. Let $p \colon \mc{Y} \to Y$ be an affinement of $Y$ of relative dimension $d(p)$,\footnote{If $Y$ is not connected, $d(p)$ is a function $\pi_0(Y) \to \Z$; we can always reduce to the case of connected $Y$, so do not dwell on this.} and choose a closed embedding $\mc{Y} \into \mathbb{A}^N$ of $\mc{Y}$ into some affine space. Let 
\[
\mathbb{A}_{\bullet}= \{ \emptyset= \mathbb{A}_{-N-1} \subset \mathbb{A}_{-N} \subset \cdots \subset \mathbb{A}_0= \mathbb{A}^N\}
\]
be a full flag of affine linear sections of $\mathbb{A}^N$, and form the Cartesian diagram
\begin{equation}\label{flagdiagram}
\xymatrix{
\mc{X}_{\bullet} \ar@/^2 pc/[rr]^{r_{\bullet}} \ar[r]^{i_{\bullet}} \ar[d] & \mc{X} \ar[r]^{p} \ar[d] & X \ar[d]^f \\
\mc{Y}_{\bullet} \ar[r] \ar[d] & \mc{Y} \ar[r]^p \ar[d] & Y \\
\mathbb{A}_{\bullet} \ar[r] & \mathbb{A}^N & \\
}
\end{equation}
We define the associated (increasing) flag filtrations
\begin{equation}
F^{\mathbb{A}_{\bullet}}_j H^*(X_{\bar{k}}, \Ql) = \ker \left(r_{-j}^* \colon H^*(X_{\bar{k}}, \Ql) \to H^*((\mc{X}_{-j})_{\bar{k}}, \Ql) \right)
\end{equation} 
and (see Remark \ref{flagexp})
\begin{equation}
F^{\mathbb{A}_{\bullet}}_j H^*_c(X_{\bar{k}}, \Ql) = \im \left( r_{!, j} \colon H^*_c((\mc{X}_j)_{\bar{k}}, \Ql) \to H^*_c(X_{\bar{k}}, \Ql) \right)
\end{equation}
Then for a general flag $\mathbb{A}_{\bullet}$, 
\begin{equation}\label{perv=flag}
\mc{P}^f_{j} H^q(X_{\bar{k}}, \Ql)= F^{\mathbb{A}_{\bullet}}_{1+d(p)-q+j} H^q(X_{\bar{k}}, \Ql),
\end{equation}
and
\begin{equation}\label{perv=flagc}
\mc{P}^f_j H^q_c(X_{\bar{k}}, \Ql)= F^{\mathbb{A}_{\bullet}}_{j-q-d(p)}H^q_c(X_{\bar{k}}, \Ql).
\end{equation}
\end{thm}
\begin{rmk}\label{flagexp}
\begin{enumerate}
\item Let us spell out the construction of the maps $r_{!, j}$. There is a canonical identification 
\[
H^k_c(X_{\bar{k}}, \Ql) \cong H^k_c(\mc{X}_{\bar{k}}, p^! \Ql) \cong H^{k+2d(p)}_c(\mc{X}_{\bar{k}}, \Ql)(d(p)),
\]
and then adjunction gives maps
\[
H^k_c((\mc{X}_{-j})_{\bar{k}}, i_{-j}^! \Ql) \to H^k_c(\mc{X}_{\bar{k}}, \Ql).
\]
As part of the definition of `general position,' we may assume the $\mc{X}_{-j}$ are smooth, so by cohomological purity these adjunction maps are identified with (Gysin) maps
\[
H^{k-2j}_c((\mc{X}_{-j})_{\bar{k}}, \Ql)(-j) \to H^k_c(\mc{X}_{\bar{k}}, \Ql).
\]
The `corestriction' maps $r_{!, -j}$ are then given by the composites
\[
H^{k+2d(p)-2j}_c((\mc{X}_{-j})_{\bar{k}}, \Ql)(d(p)-j) \to H^{k+2d(p)}_c(\mc{X}_{\bar{k}}, \Ql)(d(p)) \to H^k_c(X_{\bar{k}}, \Ql).
\]
Important for our purposes is that these are precisely the maps Poincar\'{e} dual to the pullback maps arising from the maps $\mc{X}_{-j} \to X$.
\item We need to specify what is meant by a \textit{general} flag; this will be done in \S \ref{strats}. What matters for our purposes is that there exists a Zariski open, dense subspace $\mr{Flag}^{\mr{gen}}$ inside the variety $\mr{Flag}$ of full (affine linear) flags in $\mathbb{A}^N$ such that all points $\mathbb{A}_{\bullet} \in \mr{Flag}^{\mr{gen}}(k)$ are `general.'
\item Note that the degree shift between the perverse and flag filtrations depends on the degree ($q$ above) of cohomology. We will ultimately work one degree of cohomology at a time, and all that matters for us is that \textit{some} shift of the flag filtration agrees with the perverse filtration. To extract the exact degree shift $j \mapsto 1+d(p)-q+j$, use \cite[Theorem 3.3.5, (3.8), Example 3.1.6, and (3.16)]{decataldo:pervfil}, and similarly for cohomology with compact supports.
%
\end{enumerate}
\end{rmk}
We now explain why this result can be refined $k$-rationally, so that the diagram (\ref{flagdiagram}) for which the conclusion of Theorem \ref{pervleraythm} holds can be taken to be a diagram in $\mr{Sch}/k$. In the process, we will say more explicitly what is meant by a `general' flag in the case of \textit{proper} $f \colon X \to Y$ (this case is somewhat simpler--see \cite[Remark 3.2.13]{decataldo:pervfil}--and it is all we need).
\subsection{Stratifications}\label{strats} To define `general' flags, we need to say something about stratifications. From now on we will consider a proper map $f \colon X \to Y$ of varieties over $k$ with $Y$ quasi-projective. For the purposes of Theorem \ref{pervleraythm}, we need only find a stratification $\Sigma$ of $Y$ such that $f_* \Ql$ is $\Sigma$-constructible. That is, we require a decomposition $Y= \sqcup_{\sigma \in \Sigma} Y_{\sigma}$ of $Y$ into locally closed, irreducible, non-singular varieties such that $f_* \Ql|_{Y_{\sigma}}$ is lisse for all $\sigma$. This is easily arranged; note that the strata $Y_{\sigma}$ may be irreducible but not geometrically connected. Then we can deduce:
\begin{cor}\label{kpervleray}
Let $k$ be a field of characteristic zero, and let $f \colon X \to Y$ be a proper morphism in $\mr{Sch}/k$ with $Y$ quasi-projective. Then there exists a diagram (\ref{flagdiagram}) defined over $k$ for which the conclusions (\ref{perv=flag}) and (\ref{perv=flagc}) of Theorem \ref{pervleraythm} holds.
\end{cor}
\proof
Choose as before an affinement $p \colon \mc{Y} \to Y$ and an embedding $\mc{Y} \into \mathbb{A}^N$, and let $\mr{Flag}$ denote the variety over $k$ of full affine linear flags in $\mathbb{A}^N$. Fix a stratification $\Sigma$ of $Y$ such that $f_* \Ql$ is $\Sigma$-constructible, and pull it back to a stratification $p^{-1} \Sigma$ of $\mc{Y}$. We then consider full flags 
\[
\{\mathbb{A}_{-N} \subset \cdots \subset \mathbb{A}_{-1} \subset \mathbb{A}^N\}
\]
such that $\mathbb{A}_{-1}$ intersects every stratum $\mc{Y}_{\sigma}$ transversally; and, refining each $\mathbb{A}_{-1} \cap \mc{Y}_{\sigma}$ to the disjoint union of its connected components, $\mathbb{A}_{-2}$ intersects the induced stratification of $\mc{Y} \cap \mathbb{A}_{-1}$ transversally; and so on inductively. By Bertini's theorem in exactly the form \cite[Theoreme 6.3(2)]{jouanolou:bertini}, applied inductively to each of the (smooth) strata in each $\mc{Y} \cap \mathbb{A}_{-i}$, the collection of such flags defines a Zariski open (over $k$) dense subset $\mr{Flag}^{\mr{gen}} \subset \mr{Flag}$. Since $\mr{Flag}$ is a rational variety (for instance, by Bruhat decomposition), and $k$ has characteristic zero, $\mr{Flag}^{\mr{gen}}(k)$ is non-empty. The corollary then follows by the proof of \cite[Theorem 3.3.1, Theorem 3.3.5]{decataldo:pervfil}.
\endproof
So that we can directly invoke the results of de Cataldo-Migliorini, in what follows we will make further demands on the stratification, as explained in \cite[\S 1.3.2]{decataldo-migliorini:decomproj}. For a fixed proper map $f \colon X \to Y$, we consider stratifications of $X$ and $Y$ as the disjoint unions of smooth, locally closed, irreducible (over $k$) subvarieties, such that every stratum of $X$ maps smoothly and surjectively onto a stratum of $Y$. Organizing the strata of $Y$ by dimension, we write $Y= \sqcup_{l=0}^{\dim Y} S_l$, where $S_l$ has pure dimension $l$. Each $S_l$ is a disjoint union of smooth and irreducible components of dimension $l$; these irreducible components need not be geometrically irreducible, but that does not affect our arguments. We then have Zariski open (dense) subsets $U_l= \sqcup_{m \geq l} S_m$, and we get associated closed and open immersions $\alpha_l \colon S_l \into U_l$ (closed) and $\beta_l \colon U_{l+1} \into U_l$ (open), with $U_l= S_l \sqcup U_{l+1}$. For more background on these stratifications, see \cite[\S 3.2]{decataldo-migliorini:htalgmaps}.
\subsection{Motivated perverse Leray filtration}
By Corollary \ref{kpervleray}, we have $\gal{k}$-equivariantly identified the perverse Leray filtrations on $H^*(X_{\bar{k}}, \Ql)$ and $H^*_c(X_{\bar{k}}, \Ql)$ with certain flag filtrations, given in terms of maps of $k$-varieties $X \to \mc{X}_{-j}$. With an eye toward our final application, in which case $f \colon X \to Y$ will be a resolution of singularities of $Y$, we continue to assume $f$ is proper, but also require that $X$ is non-singular and irreducible,\footnote{The irreducibility assumption is only for convenience in certain intermediate results, in which, for instance, we wish to invoke Poincar\'{e} duality without complicating the notation. Eventually, we extend component by component to the reducible case.} and we use Theorem \ref{hfunctor}, Lemma \ref{motivatedopenPD} and Corollary \ref{kpervleray} to \textit{define} the `perverse Leray filtration' on the motivated motives $H_r (h(X))$ and $H^s(W(X))$. Consider a diagram (\ref{flagdiagram}) over $k$ for which the conclusion of Theorem \ref{pervleraythm} holds. Since $h \colon \mr{Sch}/k \to K^b(\mc{M}_k)$ is a functor, we obtain a commutative diagram
\[
\xymatrix{
h(X) \ar[r]^{r_{-i}^*} \ar[dr]_{r_{-i+1}^*} & h(\mc{X}_{-i}) \ar[d] \\
& h(\mc{X}_{-i+1}) \\
}
\]
in $K^b(\mc{M}_k)$. Recall that since $\mc{M}_k$ is canonically weight-graded (it has K\"{u}nneth projectors), we can apply the composite functor $H^q_r$ given by taking cohomology $H_r$ of this diagram and projecting to the weight $q+r$-component for any $q \in \Z$, obtaining a commutative diagram in $\mc{M}_k$,
\[
\xymatrix{
H_r^q(h(X)) \ar[r]^{r_{-i}^*} \ar[rd]_{r_{-i+1}^*} & H_r^q(h(\mc{X}_{-i})) \ar[d] \\
& H_r^q(h(\mc{X}_{-i+1})).\\
}
\]
Recall that the $\ell$-adic realization of $H_r^q(h(X))$ is (by Theorem \ref{hfunctor}) isomorphic to $\Gr^W_{q+r} H^q(X_{\bar{k}}, \Ql)$; this accounts for the notation.
\begin{defn}\label{motpervleraydef}
The perverse Leray filtration of $H_r^q(h(X))$ is defined to be
\begin{equation}
\mc{P}^f_j H^q_r(h(X))= \ker \left\{ r^*_{-(1+d(p)-q+j)} \colon H^q_r(h(X)) \to H^q_r(h(\mc{X}_{-(1+d(p)-q+j)})) \right\}.
\end{equation}
The gradeds for the perverse filtration, still objects of $\mc{M}_k$, are then denoted
\[
\Gr^{\mc{P}^f}_j H^q_r(h(X))= \mc{P}^f_j H^q_r(h(X))/\mc{P}^f_{j-1} H^q_r(h(X)).
\]
\end{defn}
\begin{rmk}
Our indexing convention is somewhat different from that of de Cataldo-Migliorini (compare \cite[Definition 2.2.2]{decataldo-migliorini:htalgmaps}).
\end{rmk}
Now, we already have a definition (equation (\ref{defpervfil})) of $\mc{P}^f$ on $H^q(X_{\bar{k}}, \Ql)$ before passing to the weight-gradeds; the two versions of $\mc{P}^f$ are compatible in the following sense:
\begin{lemma}
The $\ell$-adic realization $H_{\ell}\left( \mc{P}^f_j H^q_r(h(X))\right)$ is isomorphic to $\Gr^W_{q+r} \mc{P}^f_j H^q(X_{\bar{k}}, \Ql)$, and likewise with $\Gr^{\mc{P}^f}_j$ in place of $\mc{P}^f_j$.
\end{lemma}
\proof
The $\ell$-adic realization functor is exact, and the maps on cohomology induced by the morphisms $\mc{X}_{\bullet} \to X$ are strict for the associated weight filtrations, so this follows from the choice of $\mc{X}_{\bullet}$ as in Corollary \ref{kpervleray}.
\endproof
We also need a `motivated' description of the perverse Leray filtration in compactly-supported cohomology, i.e. a filtration by sub-motives on each $H^s(W^p(X))$. Taking our cue from Remark \ref{flagexp}, we formally define a filtration on $H^q_r(h(X))^\vee$ by
\begin{equation}
\left( H^q_r(h(X))/ \mc{P}^f_j H^q_r(h(X)) \right)^\vee \subset H^q_r(h(X))^\vee,
\end{equation} 
and then invoke duality to define:
\begin{defn}
The perverse Leray filtration of $H^s(W^p(X))$ is defined to be
\begin{align}\label{motivatedlerayc}
\mc{P}^f_j H^r(W^{2 \dim X-q-r}(X))&= \left( H^q_r(h(X))/ \mc{P}^f_{-j+2\dim X-1} H^q_r(h(X)) \right)^\vee(- \dim X) \\ \nonumber
&\subseteq H^q_r(h(X))^\vee(-\dim X) \xrightarrow{\sim} H^r(W^{2 \dim X -q-r}(X)).
\end{align}
\end{defn}
Again, we check that this definition is compatible with the usual one in cohomology:
\begin{lemma}
The $\ell$-adic realization $H_{\ell} \left( \mc{P}^f_j H^r(W^{2 \dim X-q-r}(X)) \right)$ is canonically isomorphic to 
\[
\Gr^W_{2\dim X-q-r} \mc{P}^f_j H^{2\dim X-q}_c(X_{\bar{k}}, \Ql),
\]
and likewise for $\Gr^{\mc{P}^f}_j$.
\end{lemma}
\proof
By the description (Remark \ref{flagexp}) of $r_{!, -j}$ as the map Poincar\'{e} dual to $r_{-j}^*$, we see that Poincar\'{e} duality for $X$ induces a duality (here $F_{\bullet}$ denotes the flag filtrations for general flags)
\[
\left(F_j H^q(X_{\bar{k}}, \Ql) \right)^\vee \cong \left( H^{2\dim X-q}_c(X_{\bar{k}}, \Ql)/F_{-j} H^{2 \dim X-q}_c(X_{\bar{k}}, \Ql)\right) (\dim X),
\] 
i.e.
\[
\left( \mc{P}^f_{-l+2\dim X-1} H^q(X_{\bar{k}}, \Ql) \right)^\vee \cong \left( H^{2\dim X-q}_c(X_{\bar{k}}, \Ql)/\mc{P}^f_{l} H^{2 \dim X-q}_c(X_{\bar{k}}, \Ql)\right) (\dim X).
\]
The lemma follows by passing to $\Gr^W_{\bullet}$.
\endproof
By definition, we obtain the following duality in $\mc{M}_k$, a motivated analogue of \cite[\S 1.3.3(12)]{decataldo-migliorini:decomproj}:
\begin{equation}\label{pervgradedPD}
\Gr^{\mc{P}^f}_j H^r(W^{2\dim X-q-r}(X)) \times \Gr^{\mc{P}^f}_{-j+2\dim X} H^q_r(h(X)) \to \Q(-\dim X).
\end{equation}
We next check a functoriality property of these motivated perverse Leray filtrations.
\begin{lemma}\label{functoriality}
Suppose
\[
\xymatrix{
\mc{T} \ar[r] \ar[dr]_{g} & X \ar[d]^f \\
& Y \\
}
\]
is a commutative diagram in $\mr{Sch}/k$. Then the pull-back maps $H^q_r(h(X)) \to H^q_r(h(\mc{T}))$ induce morphisms (in $\mc{M}_k$)
\[
\mc{P}^f_j H^q_r(h(X)) \to \mc{P}^{g}_j H^q_r(h(\mc{T})).
\]
If $g$ factors as $\mc{T} \xrightarrow{\gamma} Z \xrightarrow{\iota} Y$ with $\iota$ a closed immersion, then the filtrations $\mc{P}^{\gamma}_{\bullet}$ and $\mc{P}^{g}_{\bullet}$ on $H^*(\mc{T}_{\bar{k}}, \Ql)$, or on $H^q_r(h(\mc{T}))$, coincide.

If $\mc{T} \to X$ is proper, then proper pull-back $H^s(W^p(X)) \to H^s(W^p(\mc{T}))$ also preserves the perverse Leray filtrations (\ref{motivatedlerayc}).
\end{lemma}
\proof
Since the $\ell$-adic realization functor on $\mc{M}_k$ is exact, it suffices to check the statement in cohomology. Here it is elementary: see for instance \cite[Remark 4.2.3]{decataldo-migliorini:htalgmaps}. For the second statement, use the fact that $\iota_*$ is exact for the perverse t-structure (so commutes with perverse truncation).
\endproof
\subsection{Motivated support decomposition}
Now we proceed as in \cite[Proposition 2.2.1]{decataldo-migliorini:decomproj} to establish a `motivated support decomposition' of the $\Gr^{\mc{P}^f}_j H^q_r(h(X))$, corresponding to the support decomposition of the perverse sheaf ${}^p H^j(f_* \Ql)$. We begin by checking that the desired support decomposition exists $k$-rationally. Continue to let $f \colon X \to Y$ be our proper map of quasi-projective varieties over $k$ with $X$ non-singular. Let $Y= \sqcup_{l=0}^{\dim Y} S_{l}$ be a stratification for $f$ as in \S \ref{strats}, with the collection of closed and open immersions
\begin{equation}
S_{l} \xrightarrow{\alpha_{l}} U_{l} \xleftarrow{\beta_{l}} U_{l+1}.
\end{equation} 
\begin{lemma}\label{supportdecomp}
In the above setting, there is a canonical isomorphism in $\mr{Perv}(Y)$
\begin{equation}\label{supporteqtn}
{}^p H^j (f_* \Ql[\dim X]) \xrightarrow{\sim} \bigoplus_{l= 0}^{\dim Y} \IC_{\overline{S_{l}}}\left( \alpha_{l}^* H^{-l} ({}^p H^j (f_* \Ql [\dim X]) \right),
\end{equation}
where the $\alpha_{l}^* H^{-l} ({}^p H^j (f_* \Ql [\dim X])$ are (geometrically semi-simple) local systems on $S_l$. Replacing $\alpha_{l}$ by the inclusion $S \xrightarrow{\iota_S} S_{l} \xrightarrow{\alpha_{l}} U_{l}$ of an irreducible ($=$ connected) component $S$ of $S_l$, we obtain the refined $k$-rational support decomposition
\begin{equation}\label{fullsupportdecomp}
{}^p H^j (f_* \Ql[\dim X]) \xrightarrow{\sim} \bigoplus_{l= 0}^{\dim Y} \bigoplus_{S \in \pi_0(S_{l})} \IC_{\overline{S}} \left(  (\alpha_{l} \circ \iota_S)^* H^{-l} ({}^p H^j (f_* \Ql [\dim X])\right).
\end{equation}
\end{lemma}
\proof
The second claim follows from the first, so we focus on establishing (\ref{supporteqtn}). This statement in $\mr{Perv}(Y_{\bar{k}})$ is a precise form--see \cite[Theorem 2.1.1c]{decataldo-migliorini:htalgmaps}--of the semi-simplicity assertion of the decomposition theorem, so it suffices to check that the map in equation (\ref{supporteqtn}) can be defined in $\mr{Perv}(Y)$. For notational simplicity, denote the perverse sheaf ${}^p H^j (f_* \Ql [\dim X])$ on $Y$ simply by $K$. We follow closely the argument of \cite[Lemma 4.1.3]{decataldo-migliorini:htalgmaps}, and, as there, the claim will follow from the following assertion: for all $l=0, \ldots, \dim Y$, there is a canonical isomorphism 
\[
K|_{U_l} \xrightarrow{\sim} \beta_{l !*}(K|_{U_{l+1}}) \oplus H^{-l}(K|_{U_l})[l].
\]
We now explain this isomorphism, which itself follows from the the corresponding geometric statement in \cite[Lemma 4.1.3, \S 6]{decataldo-migliorini:htalgmaps}. The second projection comes from the truncation triangle 
\[
\tau_{\leq -l-1} K|_{U_l} \to \tau_{\leq -l} K|_{U_l} \to H^{-l} K|_{U_l}[l] \xrightarrow{+1},
\]
whose middle term is canonically $K|_{U_l}$, and whose right-hand term is perverse (by the support conditions in the definition of perverse sheaves; see \cite[\S 4.1]{decataldo-migliorini:htalgmaps}).

To define the first projection, recall the successive truncation description of intermediate extension as
\[
\tau_{\leq -l-1} \beta_{l *}\beta_l^* K|_{U_l} \cong \beta_{l !*}(K|_{U_{l+1}}).
\]
This suggest applying $\Hom(K|_{U_l}, \cdot)$ to the truncation triangle
\[
\tau_{\leq -l-1} \beta_{l *}\beta_l^* K|_{U_l} \to \tau_{\leq -l} \beta_{l *}\beta_l^* K|_{U_l} \to H^{-l}(\beta_{l *}\beta_l^* K|_{U_l})[l] \xrightarrow{+1}.
\]
There is a canonical map 
\[
\xymatrix{
K|_{U_l} \ar@/_2 pc/[rr]^a & \tau_{\leq -l} K|_{U_l} \ar[l]_{\sim} \ar[r] & \tau_{\leq -l} \beta_{l*} \beta_l^* K|_{U_l},
}
\]
and to construct the projection $K|_{U_{l}} \to \beta_{l !*}(K|_{U_{l+1}})$, it suffices to check that the image of $a$ in 
\[
\Hom(H^{-l}(K|_{U_l}), H^{-l}(\beta_{l *}\beta_l^* K|_{U_l})) \xrightarrow{\sim} \Hom(K|_{U_l}, H^{-l}(\beta_{l *}\beta_l^* K|_{U_l})[l])
\]
(recall $\tau_{\leq -l} K|_{U_l} \xrightarrow{\sim} K|_{U_l}$) is zero. But we can check whether a map of constructible sheaves on $Y$ is zero by passing to $Y_{\bar{k}}$, so the geometric assertion (\cite[\S 6]{decataldo-migliorini:htalgmaps}) implies our corresponding arithmetic assertion.
\endproof
We have maneuvered into a position to invoke the argument of \cite[Proposition 2.2.1]{decataldo-migliorini:decomproj} to prove:
\begin{thm}\label{motivateddecomposition}
Let $f \colon X \to Y$ be a proper map of quasi-projective varieties over $k$ with $X$ non-singular. Then for each triple of integers $j$, $q$, $r$, there exists a decomposition in $\mc{M}_k$
\begin{equation}
\Gr^{\mc{P}^f}_j H^q_r(h(X)) \xrightarrow{\sim} \bigoplus_{l=0}^{\dim Y} \bigoplus_{S \in \pi_0(S_l)} \Gr^{\mc{P}^f}_{j, S} H^q_r(h(X))
\end{equation}
whose $\ell$-adic realization is the output of applying $\Gr^W_{q+r} H^q(Y_{\bar{k}}, \bullet)$ to the splitting of Lemma \ref{supportdecomp}, equation (\ref{fullsupportdecomp}).\footnote{Rather, the slight relabeling of this splitting that results from replacing $f_*\Ql[\dim X]$ in equation (\ref{fullsupportdecomp}) with $f_*\Ql$.}

Similarly for cohomology with compact supports, i.e. for the motives $\Gr^{\mc{P}^f}_j H^s(W^p(X))$.
\end{thm}
\proof
When $X$ is projective, this is established in \cite[Theorem 3.2.2]{decataldo-migliorini:decomproj}, via the argument of \cite[Proposition 2.1.1]{decataldo-migliorini:decomproj}; the impediment to that argument going through for non-projective $X$ is dealt with by our systematic use of the motivated motives $H^q_r(h(X))$ and $H^s(W^p(X))$. We do not repeat the proof, but we will remark on the key points. The argument proceeds by induction on $\dim X$; the inductive step is achieved by using \cite[Equations (13) and (14)]{decataldo-migliorini:decomproj} to define the summands $\Gr^{\mc{P}^f}_{j, S} H^q_r(h(X))$ in terms of already-defined terms for lower-dimensional $X$. Note that it is essential that we have at our disposal motives corresponding both to cohomology without supports (the $H^q_r(h(X))$) and to cohomology with compact supports (the $H^s(W^p(X))$), with their respective pull-back functorialities (Theorems \ref{hfunctor} and \ref{Wfunctor}) and the duality (Lemma \ref{motivatedopenPD} and Equation (\ref{pervgradedPD})) relating them.
\endproof
The corresponding result in \cite{decataldo-migliorini:decomproj} uses the relative hard Lefschetz theorem to obtain an absolute Hodge splitting (in the case $k= \bar{k}$) of $H^*(X_{\bar{k}}, \Ql)$ corresponding to the full splitting of $f_* \Ql$ given by the decomposition theorem, rather than as here merely the support decomposition in a particular perverse degree; a similar strengthening can be established in our context, which we now briefly sketch, although it is not needed for our primary goal, Theorem \ref{ICmotivated}. For $X$ as in the theorem, consider as usual a smooth compactification $\overline{X}$ with $\overline{X} \setminus X = \cup_{\alpha} D_{\alpha}$ equal to a union of smooth divisors with normal crossings. We may assume $\overline{X}$ is projective, and then take $\eta$ to be a hyperplane line bundle arising from some projective embedding. The required motivated version of the relative Hard Lefschetz theorem is that there are isomorphisms
\begin{equation}\label{motHL}
\cup \eta^j \colon \Gr^{\mc{P}^f}_{-j+\dim X} H^q_r(h(X)) \xrightarrow{\sim} \Gr^{\mc{P}^f}_{j+ \dim X}H^{q+2j}_r(h(X)) (j).
\end{equation}
To construct this isomorphism, note first that we can pull $\eta$ back to any of the intersections $D_{\Sigma}= \cap_{\alpha \in \Sigma} D_{\alpha}$ and obtain a morphism of complexes $h(X) \to h(X)(1)$;\footnote{Essentially by the projection formula: writing $\eta_{\Sigma}$ for the pull-back of $\eta$ to $D_{\Sigma}$, we have, for any inclusion $\iota \colon D_{\Sigma} \into D_{\Sigma'}$, $\iota_*(a \cup \eta_{\Sigma})= \iota_*(a) \cup \eta_{\Sigma'}$. That is, cup-product with $\eta$ commutes with the boundary (Gysin) maps of the complex $h(X)$.} passing to cohomology, $\eta$ induces maps
\[
\eta \colon H^q_r(h(X)) \to H^{q+2}_r(h(X))(1).
\]
The required compatibility
\[
\eta \colon \mc{P}^f_j H^q_r(h(X)) \to \mc{P}^f_{j+2} H^{q+2}_r(h(X))(1)
\]
with the perverse filtrations follows directly from Definition \ref{motpervleraydef}. We therefore have constructed the maps appearing in Equation (\ref{motHL}); that they are isomorphisms, as are the corresponding maps for each term of the support decomposition, then follows as usual from the corresponding statement in cohomology. The formalism of `hard Lefschetz triples' (\cite[\S 1.3.4]{decataldo-migliorini:decomproj}) in the abelian category $\mc{M}_k$ allows us to enhance Theorem \ref{motivateddecomposition} with the following:
\begin{cor}\label{splitting}
The choice of $\eta$ gives rise to a distinguished splitting of the motivated perverse Leray filtration, 
\[
H^q_r(h(X)) \xrightarrow{\sim} \bigoplus_j \Gr^{\mc{P}^f}_j H^q_r(h(X))
\]
\end{cor}
Of course, $\mc{M}_k$ is semi-simple, so we already knew that some splitting exists. The combination of Theorem \ref{motivateddecomposition} and Corollary \ref{splitting} can be regarded as a `motivated decomposition theorem.' Finally, we reach the motivating application:
\begin{cor}[Theorem \ref{ICmotivated} above]\label{ICmotpf}
Let $Y/k$ be any quasi-projective variety. Then there is an object $M \in \mc{M}_k$ whose $\ell$-adic realization is isomorphic as $\gal{k}$-representation to $\Gr^W_i \IH^q(Y_{\bar{k}}, \Ql)$. If $\mr{A}$ is a finite group scheme over $k$ acting on $Y$, and $e \in \Qb[\mr{A}(\bar{k})]^{\gal{k}}$, then for any embedding $\Qb \into \Qlb$ there is an object of $\mc{M}_{k, \Qb}$ whose $\Qb \into \Qlb$-realization is isomorphic as $\gal{k}$-representation to $\Gr^W_i e(\IH^q(Y_{\bar{k}}, \Qlb))$. 

The same holds for intersection cohomology with compact supports.
\end{cor}
\proof
We can assume $Y$ is irreducible. Let $f \colon X \to Y$ be a resolution of singularities; $X$ is then irreducible of dimension $\dim X$. For the motive $M$ having $\ell$-adic realization $\Gr^W_{q+r} \IH^q(Y_{\bar{k}}, \Ql)$, we can take
\[
M= \Gr^{\mc{P}^f}_{\dim X, Y^{\mr{sm}}} H^q_r(h(X)),
\]
where $Y^{\mr{sm}}$ denotes the smooth locus of $Y$. (Compare \cite[Remark 1.4.2]{decataldo-migliorini:decomproj}, noting that we have normalized the perverse filtration differently than they do.)

For the equivariant statement, take $f \colon X \to Y$ to be an $\mr{A}$-equivariant resolution of singularities.\footnote{For the existence, including our case in which $\mr{A}$ is not necessarily a discrete group scheme, see, e.g., \cite[9.1]{kollar:resolutionlecture}.} By Lemma \ref{functoriality}, each $\gamma \in \Gamma(\bar{k})$ induces an automorphism of $\Gr^{\mc{P}^f}_{\dim X} H^q_r(h(X_{\bar{k}}))$. For $e \in \Qb[\mr{A}(\bar{k})]^{\gal{k}}$, we obtain (after extending scalars to $\Qb$) an endomorphism of $\Gr^{\mc{P}^f}_{\dim X} H^q_r(h(X))$. That this endomorphism preserves the canonical sub-motive (Theorem \ref{motivateddecomposition})
\[
\Gr^{\mc{P}^f}_{\dim X, Y^{\mr{sm}}} H^q_r(h(X)) \subset \Gr^{\mc{P}^f}_{\dim X} H^q_r(h(X))
\]
is then verified by checking the corresponding statement for $\ell$-adic realizations.

The statement for compact supports follows similarly, or by now invoking Poincar\'{e} duality.
\endproof
\begin{rmk}\label{finalrmk}
\begin{enumerate}
\item The motive underlying $\Gr^W_i \IH^k(Y_{\bar{k}}, \Ql)$ is canonical in the following sense. The only ambiguity in its construction is that we may take a second resolution $f' \colon X' \to Y$ before applying the argument of Corollary \ref{ICmotpf}. But any two resolutions of singularities can be dominated by a third, and so the functoriality property of Lemma \ref{functoriality} implies that by passing through this third resolution we can deduce an isomorphism in $\mc{M}_k$: 
\[
\Gr^{\mc{P}^f}_{\dim X, Y^{\mr{sm}}} H^q_r(h(X)) \cong \Gr^{\mc{P}^{f'}}_{\dim X', Y^{\mr{sm}}} H^q_r(h(X')).
\]
\item In particular, Corollary \ref{ICmotpf} completes the proof of Theorem \ref{motivated}.
\item Let us now take $k$ to be a finite extension of $\Q_p$, with $\ell=p$. Let $Y/k$ be a projective variety--this way we avoid discussing weight filtrations, and in particular do not have to be concerned that $k$ is not finitely-generated over $\Q$--so that $\IH^q(Y_{\bar{k}}, \Q_p)$ has underlying motive $\Gr^{\mc{P}^f}_{\dim X, Y^{\mr{sm}}}(H^q(X))$ where $f \colon X \to Y$ is any resolution of singularities. By part 1 of Remark \ref{finalrmk}, we can then canonically define the intersection de Rham cohomology of $Y/k$ to be the de Rham realization (a filtered $k$-vector space) of the motive $\Gr^{\mc{P}^f}_{\dim X, Y^{\mr{sm}}}(H^q(X))$, and by general properties of $\mc{M}_k$ we obtain a $p$-adic de Rham comparison isomorphism, compatible with morphisms of motivated motives.
\item Finally, taking $k$ to be a totally real field, \cite[Corollary B]{patrikis-taylor:irr} extends from smooth projective varieties over $k$ with Hodge-regular cohomology in some degree to arbitrary projective varieties over $k$ with Hodge-regular intersection cohomology in some degree. Here we use the theorems of Gabber that $\{ \IH^q(Y_{\bar{k}}, \Ql)\}_{\ell}$ forms a weakly compatible system of pure $\gal{k}$-representations. Consequently, these compatible systems (in the regular case) are strongly compatible, and the corresponding L-functions admit meromorphic continuation to the whole complex plane, with the expected functional equation. Is it possible to construct examples of such singular varieties $Y$? Note that Yun's construction in type $G_2$ and $D_{2n}$ (the latter regarded as $\mr{SO}_{4n-1}$-valued) do give families of examples of potentially automorphic motives: this is a special case of the examples arising from Katz's theory as discussed in \cite[\S 2]{patrikis-taylor:irr}. The lifts of Yun's examples constructed in Corollary \ref{locsyslift} are no longer Hodge-Tate regular, so no further examples of potentially automorphic motives result from the constructions of the present paper.
\end{enumerate}
\end{rmk}
\bibliographystyle{amsalpha}
\bibliography{biblio.bib}

\end{document}